\documentclass[10pt]{amsart}

\usepackage{amsmath,amsthm, graphicx, amssymb}
\usepackage[all]{xy}

\usepackage{amsfonts}

\newtheorem{theorem}{Theorem}[section]
\newtheorem{proposition}[theorem]{Proposition}
\newtheorem{lemma}[theorem]{Lemma}
\newtheorem{corollary}[theorem]{Corollary}

\theoremstyle{remark}
\newtheorem{remark}[theorem]{Remark}
\newtheorem{example}[theorem]{Example}

\begin{document}

%newcommands - operations
\newcommand{\br}[1] { \left<{#1}\right>}
\newcommand{\brc}[1] { \left<{#1}\right>_c}
\newcommand{\bbr}[2] { \left<{#1}\,,\,{#2}\right>}
\newcommand{\bbrc}[2] { \left<{#1}\,,\,{#2}\right>_c}
\newcommand{\bbrx}[2] { \left<{#1}\,,\,{#2}\right>_{\mathcal{X}}}
\newcommand{\bbrxc}[2] { \left<{#1}\,,\,{#2}\right>_{\mathcal{X} , \,c  }}

\newcommand{\dbr}[1] { \left< \br{#1} \right>}
\newcommand{\dbrc}[1] { \left< \br{#1} \right>_c}
\newcommand{\dbbr}[2] { \left< \bbr{#1}{#2} \right>}
\newcommand{\dbbrc}[2] { \left< \bbr{#1}{#2} \right>_c}
% Change of notation in generalization section
\newcommand{\rbbr}[2] { \partial_{{#1}}\,\left({#2}\right)}
\newcommand{\rbbrc}[2] { \partial_{#1}\, \left({#2}\right)_c}
%\newcommand{\rbr}[1] { \left[{#1}\right]}
%\newcommand{\rbbr}[2] { \left[{#1}\,,\,{#2}\right]}
%\newcommand{\rbbrc}[2] { \left[{#1}\,,\,{#2}\right]_c}
%\newcommand{\rbrc}[1] { \left[{#1}\right]_c}

%frown
%%%%%%%%% commented for no diagram version
\newcommand{\ncstrut}{\xy (0,0.5)*{\bullet}; (2,1.2)*{\frown};
(4.6,0.5)*{\bullet~}; \endxy }
\newcommand{\cstrut}[2] { \xy (0,0)*{#1}; (2.3,0)*{\bullet}; (4.3,0.7)*{\frown}; (6.3,0)*{\bullet}; (9.4,0)*{#2~}; \endxy }

% \newcommand{\ncstrut}{\frown}
% % \newcommand{\cstrut}[2] {{#1}\frown{#2}}

%newcommands - notation
\newcommand{\g}{{\sf g}}
\newcommand{\h}{{\sf h}}
\newcommand{\Y}{{\mathcal{Y}}}
\newcommand{\X}{{\mathcal{X}}}

% im %  since \AA doesn't look nice in math mode
\newcommand{\aar}{\text{{\rm \AA}}}

% dmj
\renewcommand{\familydefault}{ppl}
\renewcommand{\thesection}{\arabic{section}}

\title{On the group-like behaviour of the Le-Murakami-Ohtsuki invariant}

\author[D.~M.~Jackson]{David~M.~Jackson$^*$}
\author[I. Moffatt]{Iain Moffatt$^\dagger$}
\author[A.H. Morales]{Alejandro Morales$^\ddagger$}

\thanks{
${\hspace{-1ex}}^*$Department of Combinatorics and Optimization,
                                                University of Waterloo, Waterloo, Ontario, Canada;  \\
${\hspace{.35cm}}$ \texttt{dmjackso@math.uwaterloo.ca}}

\thanks{
${\hspace{-1ex}}^\dagger$ Department of Applied Mathematics (KAM \& ITI), Charles University, Prague, Czech Republic; \\
${\hspace{.35cm}}$ \texttt{iain@kam.mff.cuni.cz}}

\thanks{
${\hspace{-1ex}}^\ddagger$Faculty of Mathematics,
                                                University of Waterloo, Waterloo, Ontario, Canada;  \\ %
${\hspace{.35cm}}$ \texttt{ahmorale@math.uwaterloo.ca}}

\date{16 November, 2005}

\begin{abstract}
We study the effect of Feynman integration and diagrammatic
differential operators on the structure of group-like elements
in the algebra generated by coloured vertex-oriented uni-trivalent
graphs. We provide applications of
our results to the study of the LMO invariant, a quantum invariant
of manifolds. We also indicate further situations in which our
results apply and may prove useful. The enumerative approach that we
adopt has a clarity that has enabled us to perceive a number of
generalizations.
\end{abstract}
\maketitle
\tableofcontents

%%%%%%%%%%%%%%%%%%%%%%%%%%%%%%%%%%%%%%%%%%%%%%%%%%%%%%%%%%%%%%%%%%%%%%%%%%%%
{\parskip=12pt
%%%%%%%%%%%%%%%%%%%%%%%%%%%%%%%%%%%%%%%%%%%%%%%%%%%%%%%%%%%%%%%%%%%%%%%%%%%%

% 1 -----------------------------------------------------------------------------------------------

\section[Introduction]{Introduction}

The techniques of Feynman integration and diagrammatic differential
operators play an important role in quantum topology. Roughly
speaking, these techniques involve ``gluing together'' formal power series of (coloured
uni-trivalent) graphs according to certain recipes arising from the diagrammatic formalism of perturbative Chern-Simons theory.

Feynman diagrams appear in quantum topology as equivalence classes of formal $\mathbb{Q}$-power series of coloured vertex-oriented uni-trivalent graphs.
These power series can be equipped with a commutative multiplication (given by the disjoint union) and a coproduct (the sum of all ways of ``splitting'' diagrams) and can be made into graded Hopf algebras,  denoted by $\mathcal{B}$ (see \cite{BN} for details). The primitives of these Hopf algebras are known to be power series of connected elements.
Most of the elements of $\mathcal{B}$ which are of interest in quantum topology, such as the values of the Kontsevich or Le-Murakami-Ohtsuki (hereinafter, LMO) invariants,  are known to be group-like (\cite{LMO, ohtsuki}). A well known property of graded Hopf algebras is that any group-like element may be written as the exponential of a primitive element (\cite{Abe}). This allows one to study the logarithm of quantum invariants rather than the invariants themselves.

In this paper we study the effects that Feynman integration and diagrammatic differential operators have on the structure of group-like elements in $\mathcal{B}$ through the effect on the primitives.
We provide applications of our results to the study of the LMO invariant. We also indicate further situations where our results apply and may prove useful. The enumerative approach that we have adopted has a clarity that has enabled us to perceive a number of generalizations.

Our approach is to show that these results arise naturally as a
generalization of a classical result in algebraic combinatorics.

The paper is structured as follows. In Section~\ref{background} we describe the problem and our results in a purely combinatorial language. In Section~\ref{motivation} we explain how our results relate to quantum invariants of 3-manifolds and show how to express the values of the  primitive LMO invariant in terms of those of the primitive Kontsevich invariant.
The enumerative preliminaries are given briefly in Section~\ref{S:EP}, and the details of the labelling
process are given in Section~\ref{S:LG}.
Section~\ref{S:GSEGS} deals with the graph-subgraph series for appropriately weighted graphs.
The proof of the main theorem appears in Section~\ref{S:PMT}.
In Section~\ref{generalization} we explain how our results can be generalized and applied to diagrammatic differential operators. We conclude by applying our results to find  closed formulae for the primitive LMO invariant of certain 3-manifolds. These appear in Section~\ref{S:EX}.

% 2 -----------------------------------------------------------------------------------------------

\section[The Combinatorial Problem]{The Combinatorial Problem}\label{background}

A {\em $\Y$-coloured uni-trivalent diagram} is a graph $g$ made of
undirected edges with two types of vertices {\bf i)}~trivalent vertices equipped with a cyclic ordering of its incident edges
and {\bf ii)}~univalent vertices with colours
assigned from a finite set $\Y = \{y_1, y_2, \ldots \}$, where
vertices are to be regarded as mutually distinguishable. If
$\Y=\emptyset$ then the graph is {\em trivalent}.

For example, below are two graphs, $g_0$ and $h_0$, each with two
components. The graph $g_0$ contains both trivalent and univalent
vertices with colour set $\Y=\{y_1, y_2\}$, while the graph $h_0$
has only trivalent vertices, so $\Y=\emptyset$.

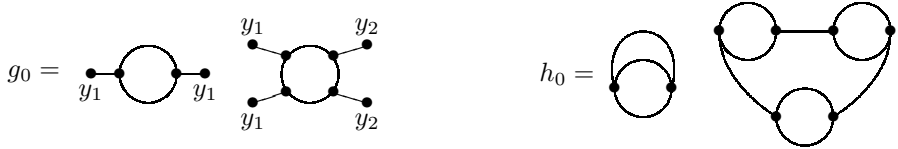
\[
\xy
% shrink
0;/r.30pc/:
%graph 1
(-6,0)*{g_0 = }; (0,-2)*{y_1}; (12,-2)*{y_1}; (0,0)*{}="A";
(3,0)*{}="B"; (9,0)*{}="C"; (12,0)*{}="D"; (0,0)*{\bullet};
(3,0)*{\bullet}; (9,0)*{\bullet}; (12,0)*{\bullet}; "A"; "B"
**\dir{-}; "C"; "D" **\dir{-}; (6,0)*\xycircle(3,3){-}="f"; "f";
%graph2
(17,5)*{y_1};
(17,-5)*{y_1};
(29,5)*{y_2};
(29,-5)*{y_2};
(17,3)*{}="A"; (20.5,1.9)*{}="B"; (17,-3)*{}="C"; (20.5,-1.9)*{}="D";
(29,3)*{}="E"; (25.5,1.9)*{}="F"; (29,-3)*{}="G"; (25.5,-1.9)*{}="H";
(17,3)*{\bullet}; (20.5,1.9)*{\bullet}; (17,-3)*{\bullet}; (20.5,-1.9)*{\bullet};
(29,3)*{\bullet}; (25.5,1.9)*{\bullet}; (29,-3)*{\bullet}; (25.5,-1.9)*{\bullet};
"A"; "B" **\dir{-};
"C"; "D" **\dir{-};
"E"; "F" **\dir{-};
"G"; "H" **\dir{-};
(23,0)*\xycircle(3,3){-}="f";
"f";
%\graph3
(50,0)*{h_0 = }; (58,0)*{\xy
    (-6,0)*{}="A"; (0,0)*{}="B";
    (-6,0)*{\bullet}; (0,0)*{\bullet};
    (-3,0)*\xycircle(3,3){-}="f1";
    "A";"B" **\crv{(-7,4) & (-3,8) & (1,4)};
    "f1";
    \endxy};
%graph4
(75,0)*{ \xy (2,0)*\xycircle(3,3){-}; (14,0)*\xycircle(3,3){-};
(8,-9)*\xycircle(3,3){-};  (5,0)*{}="X1"; (11,0)*{}="X2";
 "X1"*{\bullet}; "X2"*{\bullet};
"X1";"X2" **\dir{-};
%"X1";"X2" **\crv{(6,4)};
(-1,0)*{}="X1"; (5,-9)*{}="X2";
 "X1"*{\bullet}; "X2"*{\bullet};
"X1";"X2" **\crv{(-1,-4) & (2,-7) }; (17,0)*{}="X1";
(11,-9)*{}="X2";
 "X1"*{\bullet}; "X2"*{\bullet};
"X1";"X2" **\crv{(17,-4) & (14,-7) };
\endxy
};
\endxy
\]

Note for example, that the two vertices with the colour $y_1$ on each
component of $g_0$ are distinguishable.

Let $\mathcal{D}(\Y)$ be the algebra of formal power series with
$\Y$-coloured uni-trivalent graphs as indeterminates and
coefficients in $\mathbb{Q}$. The empty graph is allowed, and is
denoted by $1$ in the algebra. Commutative multiplication is given
by disjoint union.
Motivated by the Hopf algebra $\mathcal{B}$, we say that an element of $\mathcal{D}(\Y)$ is
{\em primitive} if it is a sum of connected graphs and {\em group-like} if it is the exponential of a primitive.
$\mathcal{D}_s(\Y)$ is the subalgebra of $\mathcal{D}(\Y)$
containing no components of the form $\ncstrut$ (called  {\em struts})
with the same colour at its vertices.

%%%%%new%%%%

Let $g$ and $h$ be  $\Y$-coloured uni-trivalent graphs. We define a  bilinear operator
$\bbr{\cdot}{\cdot}:\mathcal{D}(\Y)\otimes \mathcal{D}_s(\Y)\rightarrow \mathcal{D}(\emptyset)$,
 such that $\bbr{g}{h}$ is the sum
of all ways of identifying each of the $y$-coloured vertices in $g$
with each of the $y$-coloured vertices in $h$, for all colours $y\in
\Y$. If the numbers of $y$-coloured univalent vertices in the two
graphs do not match for some $y \in \Y$, the sum is zero. Coloured
univalent vertices that are not to be joined under the pairing are
indicated by {\em open vertices} (this will be used in Section
\ref{generalization}). We extend this bilinearly to all of
$\mathcal{D}(\Y) \otimes \mathcal{D}_s(\Y)$. Once the components of
$g$ and $h$ are joined by $\bbr{\cdot}{\cdot}$, we consider them as
{\em subgraphs} of the resulting trivalent graph. For all $D_1, D_2
\in \mathcal{D}(\Y) \otimes \mathcal{D}_s(\Y)$, we denote the
primitive part of $\bbr{D_1}{D_2}$ by $\bbrc{D_1}{D_2}$. The
following is an example of this operator.

% presentation graph 1
\xy
% shrink
0;/r.35pc/: (-4,0)*{\text{Let } D_1  =  }; (7,0)*{\xy (0,0)*{}="A";
(2,0)*{}="B";
          (6,0)*{}="C"; (8,0)*{}="D";
      "A"*{\bullet}; "B"*{\bullet}; "C"*{\bullet}; "D"*{\bullet};
      "A"; "B" **\dir{-}; "C"; "D" **\crv{~*=<0.5pt>{.}(7,0)};
      (4,0)*\xycircle(2,2){-}="f";
    "f";
    (1,-2)*{y_1};
    (8,-2)*{y_2};
    \endxy
    };
(23,0)*{\xy
    (0,0)*{}="A"; (2,0)*{}="B"; (6,0)*{}="C"; (10,0)*{}="D";
    (14,0)*{}="E"; (16,0)*{}="F"; "A"; "B" **\dir{-}; "C"; "D"
    **\dir{-}; "E"; "F" **\dir{-}; (4,0)*\xycircle(2,2){-}="f1";
    (12,0)*\xycircle(2,2){-}="f2";  "f1"; "f2";
    (1,-2)*{y_1};
    (16,-2)*{y_2};
    "A"*{\bullet}; "B"*{\bullet};"C"*{\bullet}; "D"*{\bullet};
     "E"*{\bullet}; "F"*{\bullet};
    \endxy
    };
(39,0)*{ \text{ and } D_2 = }; (50,0)*{\xy (0,0)*{}="A";
    (2,0)*{}="B";
    (6,0)*{}="C"; (8,0)*{}="D";
    "A"*{\bullet}; "B"*{\bullet}; "C"*{\bullet}; "D"*{\bullet};
    "A"; "B" **\dir{-}; "C"; "D" **\dir{-};
    (4,0)*\xycircle(2,2){-}="f"; "f";
    (1,-2)*{y_1};
    (8,-2)*{y_2};
    \endxy
    };
(62,0)*{\xy (0,0)*{}="A"; (2,0)*{}="B";
    (6,0)*{}="C"; (8,0)*{}="D"; "A"; "B" **\dir{-};
     "A"*{\bullet}; "B"*{\bullet}; "C"*{\bullet}; "D"*{\bullet};
    "C"; "D" **\dir{-};
    (4,0)*\xycircle(2,2){-}="f"; "f";
    (1,-2)*{y_1};
    (8,-2)*{y_2};
    \endxy
    };
(67,0)*{,};
%line 2
(0,-12)*{\text{then } \bbr{D_1}{D_2} = }; (15,-12)*{\xy
     (-8,-3)*{2};
    (-5,0)*{}="A"; (-1,0)*{}="B";
    (-5,-6)*{}="C"; (-1,-6)*{}="D";
    (-5,0)*{\bullet}; (-1,0)*{\bullet};
    (-5,-6)*{\bullet}; (-1,-6)*{\bullet};
    (-3,0)*\xycircle(2,2){-}="f1";
    (-3,-6)*\xycircle(2,2){-}="f2";
    "A";"C" **\crv{(-7,-3)};
    "B";"D" **\crv{(1,-3)};
    "f1"; "f2";
    \endxy
       };
(26,-12)*{ \xy (2,0)*\xycircle(2,2){-}; (10,0)*\xycircle(2,2){-};
(6,-7)*\xycircle(2,2){-}; (4,0)*{}="X1"; (8,0)*{}="X2";
 "X1"*{\bullet}; "X2"*{\bullet};
"X1";"X2" **\dir{-};
%"X1";"X2" **\crv{(6,4)};
(0,0)*{}="X1"; (4,-7)*{}="X2";
 "X1"*{\bullet}; "X2"*{\bullet};
"X1";"X2" **\crv{(1,-3) & (2,-5) }; (12,0)*{}="X1"; (8,-7)*{}="X2";
 "X1"*{\bullet}; "X2"*{\bullet};
"X1";"X2" **\crv{(11,-3) & (10,-5) };
\endxy
}; (35,-12)*{+}; (49,-13)*{ \xy (-4,-3.5)*{2};
(2,0)*\xycircle(2,2){-}; (10,0)*\xycircle(2,2){-};
(6,-9)*\xycircle(2,2){-}; (0,-6)*\xycircle(2,2){-};
(12,-6)*\xycircle(2,2){-}; (4,0)*{}="X1"; (8,0)*{}="X2";
 "X1"*{\bullet}; "X2"*{\bullet};
"X1";"X2" **\dir{-}; (0,0)*{}="X1"; (-2,-6)*{}="X2";
 "X1"*{\bullet}; "X2"*{\bullet};
"X1";"X2" **\crv{(-3,-3)}; (12,0)*{}="X1"; (14,-6)*{}="X2";
 "X1"*{\bullet}; "X2"*{\bullet};
"X1";"X2" **\crv{(15,-3)}; (12,-8)*{}="X1"; (8,-9)*{}="X2";
 "X1"*{\bullet}; "X2"*{\bullet};
"X1";"X2" **\crv{(10,-9)}; (0,-8)*{}="X1"; (4,-9)*{}="X2";
"X1"*{\bullet}; "X2"*{\bullet}; "X1";"X2"
**\crv{(2,-9)};
\endxy
 }; (60,-14)*{.};
%line 3
(12,-26)*{\text{The primitive part of this is } \bbrc{D_1}{D_2} = };
(46,-27)*{ \xy (-4,-3.5)*{2}; (2,0)*\xycircle(2,2){-};
(10,0)*\xycircle(2,2){-}; (6,-9)*\xycircle(2,2){-};
(0,-6)*\xycircle(2,2){-}; (12,-6)*\xycircle(2,2){-}; (4,0)*{}="X1";
(8,0)*{}="X2";
 "X1"*{\bullet}; "X2"*{\bullet};
"X1";"X2" **\dir{-}; (0,0)*{}="X1"; (-2,-6)*{}="X2";
 "X1"*{\bullet}; "X2"*{\bullet};
"X1";"X2" **\crv{(-3,-3)}; (12,0)*{}="X1"; (14,-6)*{}="X2";
 "X1"*{\bullet}; "X2"*{\bullet};
"X1";"X2" **\crv{(15,-3)}; (12,-8)*{}="X1"; (8,-9)*{}="X2";
 "X1"*{\bullet}; "X2"*{\bullet};
"X1";"X2" **\crv{(10,-9)}; (0,-8)*{}="X1"; (4,-9)*{}="X2";
"X1"*{\bullet}; "X2"*{\bullet}; "X1";"X2"
**\crv{(2,-9)};
\endxy
}; (57,-28)*{.};
\endxy

\noindent Note that it is necessarily linear in the indeterminates.
The following is the main theorem  of the paper. It determines the
primitive structure of $\bbr{D_1}{D_2}$ and its relation to the LMO
invariant will be discussed in the next section. We will also use
the equation to find formulae for the primitive LMO invariant of
certain manifolds in Section~\ref{S:EX}.
\begin{theorem}\label{thm2}
Let $B \in \mathcal{D}(\Y)$ and $C \in \mathcal{D}_s(\Y)$, be
primitive. Then
\begin{equation*}
\bbr{\exp B}{\exp C} = \exp\bbrc{\exp B}{\exp C}.
\end{equation*}
\end{theorem}

Before we  discuss  the relation of this theorem with quantum topology we  highlight two important corollaries and a generalization of the theorem. We note that another generalization, which was motivated by the theory of diagrammatic differential operators, will be  discussed in section~\ref{generalization}.

The first special case of this theorem occurs when $B$ consists entirely of struts.  This was the motivating example for this paper and we will see in the following section how it relates to the LMO invariant.

\begin{corollary} \label{cor1}
Let $C \in \mathcal{D}_s(\Y)$ be strutless and primitive. Let $\Y=
\{ y_1,y_2, \cdots \}$, $r_{i,j}(C)\in\mathbb{Q}$ and let
$S=\exp{(\sum_{i\geq j\geq 1} r_{i,j}(C) ~\cstrut{y_i}{y_j})}.$ Then
\begin{equation*}
\bbr{S}{\exp C} = \exp\bbrc{S}{\exp C}.
\end{equation*}
\end{corollary}

A linear operator
$\br{\cdot}:\mathcal{D}_s(\Y) \rightarrow \mathcal{D}(\emptyset)$
also arises in quantum topology (eg. \cite{BGRT:AIII}).
For some $g \in \mathcal{D}_s(\Y)$, $\br{g}$
 is defined as   the sum of all ways
of identifying pairwise all the $y$-coloured univalent vertices of
$g$ for all $y \in \Y$. If $g$ has an odd number of $y$-coloured
univalent vertices for some colour $y \in \Y$, then $\br{g}=0$. As with our previous operator, let $\brc{D}$ denote the primitive part of $\br{D}$.
\begin{corollary} \label{thm1}
Let $C\in \mathcal{D}_s(\Y)$ be strutless and primitive. Then
\begin{equation*}
\br{\exp C} = \exp\brc{\exp C}.
\end{equation*}
\end{corollary}
\begin{proof}
This follows from the observation $\br{D} =  \langle \exp(\frac{1}{2} \sum_{y\in \Y} \cstrut{ y}{y}),D
\rangle$ and the theorem.
\end{proof}
We note that Garoufalidis made a conjecture of the above form
(\cite{Ga}).

The following minor generalization of our theorem will allow us to extend our results to the LMO invariant of links in manifolds.
Let $g\in \mathcal{D}_s(\Y)$   and $h\in \mathcal{D}(\Y)$. Further suppose that $\X \subset \Y$. Then the definition of $\bbr{\cdot}{\cdot}$ can easily be extended to allow  the situation where we only glue together the univalent vertices of $g$ and $h$ whose colours are in $\X$. More precisely, let
$\bbrx{\cdot}{\cdot}:\mathcal{D}(\Y)\otimes \mathcal{D}_s(\Y)\rightarrow \mathcal{D}(\Y - \X)$,
be the operator such that $\bbrx{g}{h}$ is the sum
of all ways of identifying each of the $x$-coloured vertices in $g$
with each of the $x$-coloured vertices in $h$, for all colours $x\in
\X$. If the numbers of $x$-coloured univalent vertices in the two graphs do not match
for some $x \in \X$, the sum is zero. We extend this bilinearly to
all of $\mathcal{D}(\Y) \otimes \mathcal{D}_s(\Y)$.
Again $\bbrxc{\cdot}{\cdot}$ denotes the primitive part of $\bbrx{\cdot}{\cdot}$.
The following generalizes Theorem~\ref{thm2}.
\begin{theorem}\label{thm3}
Let $\X \subset \Y$, $B \in \mathcal{D}(\Y)$ and $C \in \mathcal{D}_s(\Y)$, be
primitive. Then
\begin{equation*}
\bbrx{\exp B}{\exp C} = \exp\bbrxc{\exp B}{\exp C}.
\end{equation*}
\end{theorem}
A proof of this statement can be obtained by a simple modification of the proof of Theorems~\ref{thm2}   and is therefore excluded.

\begin{remark} \label{rem1}
The requirement that $C \in \mathcal{D}_s(\Y)$ ensures that the
calculation of $\bbr{\exp B}{\exp C}$ and its specializations is finite for any given number of vertices.
\end{remark}

\begin{remark} \label{rem2}
In fact, the proofs of Theorems \ref{thm2},  \ref{thm3} and its generalization in Section~\ref{generalization} hold even if we remove the
requirement that $B$ and $C$ are uni-trivalent, but since we take
our motivation from quantum topology we do not work in this
generality here.
\end{remark}

% 3 -----------------------------------------------------------------------------------------------

\section[Motivation from Quantum Topology]{Motivation from Quantum Topology}\label{motivation}
Our primary motivation for this study comes from the theory of the LMO invariant, $Z^{LMO}$, introduced in \cite{LMO}. This  is a universal perturbative invariant of rational homology 3-spheres (see \cite{Oh2, BGRT:AII, ohtsuki} ), and a universal finite type invariant of integral homology 3-spheres (\cite{Le}). (Recall that for a ring $R$, a {\em $R$-homology sphere} is a 3-manifold $M$ such that $H_* (M;R ) =  H_* (S^3;R)$.)

The LMO invariant was first derived by considering the behavior of the Kontsevich integral of framed links under the two Kirby moves.
A few years later in \cite{BGRT:AI, BGRT:AII, BGRT:AIII},  the diagrammization of a physical argument led to a reformulation of the LMO invariant. This approach uses the notions of ``diagrammatic integration'' and the construction  is sometimes known as the {\em \AA rhus integral}.   Our motivation comes from this formulation of the LMO invariant. We  sketch the construction of this invariant.

Two of the fundamental spaces in quantum
topology are the coalgebras
$\mathcal{A}$ of formal $\mathbb{Q}$ power series of uni-trivalent
graphs with oriented trivalent vertices and whose uni-valent
vertices lie on an oriented compact coloured 1-manifold modulo certain relations, and the
coalgebra $\mathcal{B}$  which is the quotient space of $\mathcal{D}(\Y)$  generated by some relations.  We need not worry about the exact form of these relations here.
When there are fewer than two copies of $S^1$ in the 1-manifold, a (in
general non-commutative) multiplication on $\mathcal{A}$ is given by
connect summing copies of $S^1$ and ``stacking'' copies of the interval so that the colours match (this operation corresponds to the usual composition of tangles). The multiplication on $\mathcal{B}$ is
given by disjoint union. In fact, in such an instance, one may  make
 $\mathcal{A}$ and $\mathcal{B}$ into Hopf algebras.
We will also need to make use of the coalgebra isomorphism  $\chi : \mathcal{B} \rightarrow \mathcal{A}$, which is defined to be the average of all ways of placing all of the univalent vertices of an element of $\mathcal{B}$ onto the 1-manifold of $\mathcal{A}$ such that the $y$-coloured vertices lie on the corresponding component of the 1-manifold, for all $y \in \Y$. This is known as the {\em Poincar\'{e}-Birkhoff-Witt (PBW) isomorphism} as it is the diagrammization of the map from the theory of Lie algebras.
Details of these algebras and the theory of finite-type invariants can be found in~\cite{BN}.

Now let a framed link $L$ represent a rational homology sphere $M$ by surgery. Further suppose that the components of $L$ are in correspondence with a set $\Y$.
 It is a well known fact that  the image of the Kontsevich integral of a framed link under the inverse $\sigma$ of the PBW isomorphism  may
be written in the form $ \sigma \check{Z}(L) = \exp ( \sum_{x,y \in \Y} \frac{1}{2}
l_{xy}  \; \cstrut{x}{y} + C )$, where $C \in \mathcal{B}_s(Y)$ is
primitive and strutless; $\check{Z}$ is the Kontsevich integral normalized as in \cite{LMMO} and  $l_{xy}$ denotes the  linking number.
Because of this we can separate the struts and, using the definition
of the inner product in Section~\ref{background}, define  $Z^{LMO}_0 (L)$ by
\begin{equation} \label{arhus}
Z^{LMO}_0 (L) := \left\langle   \exp  ( \sum_{x,y \in \Y} - \frac{1}{2}
l^{xy} \; \cstrut{x}{y} ) , \exp (C) \right\rangle \in
\mathcal{B}(\emptyset),\end{equation}
 where
$(l^{xy})$ is the inverse of the linking matrix $(l_{xy})$ (note that since we restrict to rational homology spheres the linking  matrix is
non-singular). This procedure of gluing  together the terms of the
Kontsevich integral $\sigma Z(L)$ is known as {\em formal Gaussian
integration} (so called as it is the diagrammization of a perturbed Gaussian integral). The function $Z^{LMO}_0$ we have just defined is invariant
under only a handle-slide.   To make it invariant under stabilization and therefore into an invariant of 3-manifolds
requires the usual trick of normalizing by eigenvalues, and the {\em LMO invariant} is defined by
\[
Z^{LMO} (M) = Z^{LMO}_0(U_{+})^{-e_+} Z^{LMO}_0(U_{-})^{-e_-} Z^{LMO}_0 (L),
\]
where $U_{\pm}$ is the $\pm 1$ framed unknot and $e_{\pm}$ is the number of $\pm$ve eigenvalues of the linking matrix.
One immediately notices that the definition of $Z^{LMO}_0$ is of the
form of Corollary~\ref{cor1} and we obtain the following.
\begin{proposition} \label{aagplike}
The value of the LMO invariant $Z^{LMO} (M)$ is group-like in $\mathcal{B} (\emptyset)$, that is $Z^{LMO} (M) = \exp (C)$, where $C$ is primitive.
\end{proposition}
We will discuss this group-like property of the LMO invariant further in Section~\ref{Examples}.

The LMO invariant can be easily extended to tangles and links in manifolds.
A framed tangle  in a rational homology sphere can be represented by a framed tangle $T \subset S^3$ through surgery. Some of the components of this will be distinguished as {\em surgery components}, \emph{ie}. surgery along these components recovers the original manifold and the remaining components correspond to components of the tangle. Suppose that $T$ is $\Y$-coloured and the surgery components are $\X$-coloured. Then one can construct the LMO invariant of tangles in rational homology spheres in a similar way to the construction outlined above except using the operation $\bbrx{\cdot}{\cdot}$ in place of $\bbr{\cdot}{\cdot}$ and restricting the linking matrix to the surgery components of $T$. See \cite{BGRT:AII, Mof} for details. Theorem~\ref{thm3} then gives:
\begin{proposition}
The value of the LMO invariant of tangles in rational homology spheres is group-like.
\end{proposition}
As an example, this property was used  in \cite{Mof} to relate the tree part of the LMO invariant of links in integral homology spheres to Milnor's $\mu$-invariants, which are classical link invariants defined through the fundamental group of the link complement.

The group-like structure is a fundamental property of the Kontsevich and LMO invariants. The fact that the Kontsevich integral is group-like is well known. The proof of the  group-like property of the LMO invariant using Le, Murakami and Ohtshuki's construction has a very different flavor than that presented above. It is reduced to the problem of showing that  a certain diagram of  algebras is commutative. This proof and that for the Kontsevich integral can be found in \cite{LMO}.
One advantage of our proofs for the group-like property of the  LMO invariant is that it expresses their values in terms of the values of the primitive Kontsevich integral in a particularly neat way.

We note that Theorem~\ref{thm2} is much more general than what is needed for the applications above. In Section~\ref{Examples} we shall need  the theorem in its more general form.

Before we return to the combinatorics, we  briefly describe one more situation where our formula applies and may prove useful. For brevity, in this paragraph we will assume that the 1-manifold in $\mathcal{A}$ is connected and the colouring set of  $\mathcal{B}$ has exactly one element.
It was noted earlier that  the Poincar\'{e}-Birkhoff-Witt isomorphism gives a vectorspace isomorphism
$\chi : \mathcal{B} \rightarrow \mathcal{A}$. This is not an algebra isomorphism, however the existence of a product on $\mathcal{B}$ which would make this map into an algebra isomorphism is immediate.  The multiplication was calculated explicitly in \cite{BGRT:AII} (although its existence had been used in several places before) and is given in terms of gluing rooted forests. Explicitly, the multiplication is defined by the formula
\[\sigma ( \chi (D_1) \cdot \chi (D_2)) = \langle \exp \Lambda , D_1 \cdot
D_2 \rangle,\]
where $\Lambda$ is the Baker-Campbel-Hausdorf formula (this measures the failure of the identity $e^{x+y} = e^x e^y$ in a Lie algebra) (\cite{Ja}) written as rooted
trees. We refer the reader to \cite{BGRT:AII} for details. Again notice
that Theorem~\ref{thm2} applies to this situation and  gives a formula for the primitive values.

%%%%%%%%%%%%%%%
%%%%%%%%%% IAIN HAS CHANGED ABOVE: March 28th
%%%%%%%%%
%%%%%%%%%

% 4 -----------------------------------------------------------------------------------------------

\section[Enumerative Preliminaries]{Enumerative Preliminaries}\label{S:EP}

%Enumerative Preliminaries

For completeness we include some familiar elementary results. The
reader who is familiar with this can pass over this section. A more
detailed account is included in \cite{GJ}.

\subsection[Labelled structures]{Labelled structures}

The following notation will be used: $\mathcal{N}_{m} = \{
1,2,\ldots, m \}$ and for $p\geq 0$, $\mathcal{N}_{\geq p} = \{ p,
p+1, \ldots \}$. Let $f(x) = \sum_{n\geq 0} a_n x^{n}$ be a  formal
power series with coefficients in $\mathbb{Q}$. For $n\geq 0$ the
mapping $[x^n]\colon \mathbb{Q}[[x]] \rightarrow \mathbb{Q}: f(x)
\mapsto a_n$ is the {\em coefficient operator}, it is linear.

Let $\mathfrak{A}$ be a set of {\em combinatorial structures}, and
let $\varepsilon$ denote the null structure in this set. Let
$\omega(A)$ be the {\em weight} of $\mathfrak{A}$, that is, a
function $\omega:\mathfrak{A} \rightarrow \mathcal{N}_{\geq 0}$. Let
$a(n) = | \{ A \in \mathfrak{A} ~:~ \omega(A)=n \}|$. We would like
to determine this number for all $n\geq 0$. This is an {\em
enumerative problem} which we denote by $(\mathfrak{A},\omega)$. A
weight function can be {\em refined} to record more information
about a combinatorial structure by tensoring the univariate weight
functions for each item of information. So if $\omega_1, \ldots,
\omega_r: \mathfrak{A} \rightarrow \mathcal{N}_{\geq 0}$ are weight
functions, then
$$
\omega_1 \otimes \cdots \otimes \omega_r : \mathfrak{A} \rightarrow
\mathcal{N}_{\geq 0} \times \cdots \times \mathcal{N}_{\geq 0}: A
\mapsto (\omega_1(A), \ldots , \omega_r(A)).
$$
We shall need labelled structures. If $A$ is a combinatorial
structure with generic subobjects ${\bf s}$ ({\em ${\bf
s}$-subobjects}), an {\em ${\bf s}$-labelling} of $A$ is $A$
together with an assignment of the numbers $1$ to $n$ to its ${\bf
s}$-subobjects. For example, in the permutation with $1$-line
presentation $(2,1,3)$, the {\sf s}-subobjects are the positions
onto which labels may be placed. The labels on positions $1,2,3$ are
$2,1,3$. Throughout, the term ``label" and ``labelling" will be used
exclusively in this strict combinatorial sense.

\subsection[Elementary counting lemmas]{Elementary counting
lemmas}\label{lemmas}

The {\em ordinary generating series} for the enumerative
problem$(\mathfrak{A},\omega)$ in the indeterminate $x$ is
$\sum_{A\in \mathfrak{A}} x^{\omega(A)}$ and is denoted by
$[(\mathfrak{A},\omega)]_o$. Thus $a(n)$ is given by
$[x^n]~[(\mathfrak{A},\omega)]_{o}$. The {\em exponential generating
series} for $(\mathfrak{A},\omega)$ is $\sum_{A \in \mathfrak{A}}
x^{\omega(A)}/\omega(A)!$, and is denoted by
$[(\mathfrak{A},\omega)]_{e}$. Thus $a(n)$ is given by $[x^{n}/n!]
[(\mathfrak{A},\omega)]_{e}$. This is used when the elements $A$ of
$\mathfrak{A}$ are labelled structures. With additional
indeterminates, the corresponding generating series is multivariate
with a weight function of the form  $\omega_1 \otimes \cdots \otimes
\omega_r$. Moreover, a multivariate generating series may have some
indeterminates that mark some information {\em ordinarily} and
others that mark some information {\em exponentially}.

We give a brief account of the properties of ordinary and exponential series in terms
 of elementary operations on sets. These operations are the Cartesian product,
 the $\star$-product and composition with respect to these. The operations arise
 very naturally from a combinatorial point of view in the decomposition of sets of structures into their constituents.

With unlabelled structures we use the Cartesian product of sets.
With labelled structures, another set product is required, namely
one that distributes labels in all possible ways. Let $\mathfrak{A}$
and $\mathfrak{B}$ be sets of labelled combinatorial structures. Let
$A \in \mathfrak{A}$ have $k$ ${\bf s}$-subobjects and $B\in
\mathfrak{B}$ have $n-k$ ${\bf t}$-subobjects. Let
$\alpha=\{\alpha_1, \alpha_2,\ldots, \alpha_k\}$ be a set of
distinct positive integers. Without loss of generality, $1\leq \alpha_1 <
\ldots < \alpha_k \leq n$. Then $(A)_{\alpha}$ denotes the structure
obtained from $A$ by replacing canonically the label $i$ with
$\alpha_i$ for $i=1,\ldots, k$. Let $\beta = \mathcal{N}_n -
\alpha$. Then $((A)_\alpha,(B)_{\beta})$ is a labelled structure
with $n$ subobjects. Let $\mathfrak{A} \star \mathfrak{B}$ denote
the set of all $((A)_{\alpha}, (B)_{\beta})$ for all choices of
$\alpha$ as a subset of $\mathcal{N}_n$, for all pairs $(A,B)\in
\mathfrak{A}\times \mathfrak{B}$ and where $n$ takes all values
greater than or equal to one. It is understood that the contribution
of $((A)_\alpha, (B)_\beta)$ is non-null only when $\alpha$ and $\beta$ are disjoint with $|\alpha|$ and
$|\beta|$ equal to the number of ${\bf s}$-subobjects and ${\bf
t}$-subobjects of $A$ and $B$, respectively (in particular $\alpha \uplus \beta = \mathcal{N}_n$, where $\uplus$ denotes disjoint union). We have the following
familiar counting lemma.

\begin{lemma}[The $\star$-Product Lemma] \label{ProdLemma}

Let $\mathfrak{A}$ and $\mathfrak{B}$ be sets of labelled
combinatorial structures with generic ${\bf s}$-subobjects and ${\bf
t}$-subobjects. Let $\omega_{\bf s}(A)$ be the number of ${\bf
s}$-subobjects of $A\in\mathfrak{A}$ and $\omega_{\bf t}(B)$ be the
number of ${\bf t}$-subobjects in $B \in \mathfrak{B}$. If ${\bf r}$
denotes the generic subobject of $\mathfrak{A}\star \mathfrak{B}$,
that is either an ${\bf s}$-subobject or a ${\bf t}$-subobject, then
$$
[(\mathfrak{A}\star \mathfrak{B}, \omega_{\bf r})]_{e} =
[(\mathfrak{A}, \omega_{\bf s})]_{e}\cdot[(\mathfrak{B}, \omega_{\bf
t})]_{e}.
$$
\end{lemma}

We shall require two auxiliary sets:  $\mathfrak{O}({\bf o}) =\{
(1,2,\ldots,k): k = 0,1,2, \ldots\}$ is the set of all {\em
canonical ordered sets} where ${\bf o}$ is the generic subobject.
Similarly, $\mathfrak{U}({\bf u}) = \{\{1,\ldots, k\} : k =
0,1,2,\ldots \}$ is the set of all {\em canonical unordered sets}
where ${\bf u}$ is the generic subobject. Trivially, we have the
generating series:
\begin{equation}
\label{elsets} [(\mathfrak{O}, \omega_{\bf o})]_{e}(x) =
\frac{1}{1-x} \text{\quad and \quad} [(\mathfrak{U},\omega_{\bf
u})]_{e}(x) = \exp x.
\end{equation}
We denote by $\mathfrak{A} \circledast \mathfrak{B}$ the composition
of sets $\mathfrak{A}$ and $\mathfrak{B}$ with respect to the
$\star$-product, where each generic subobject of an element $A \in
\mathfrak{A}$ is replaced by an element $B \in \mathfrak{B}$ in a
unique way. We have immediately the following lemma,

\begin{lemma}[The Composition Lemma] \label{CompLemma}
Let $\mathfrak{A}({\bf s})$ and $\mathfrak{B}({\bf t})$ be sets of
structures with generic ${\bf s}$ and ${\bf t}$ subobjects,
respectively. If $\omega_{\bf s}$ and $\omega_{\bf t}$ are the
weight functions of $\mathfrak{A}$ and $\mathfrak{B}$, then
$$
[(\mathfrak{A}\circledast \mathfrak{B}, \omega_{\bf t})]_{e} =
[(\mathfrak{A},\omega_{\bf s})]_{e} \circ [(\mathfrak{B},
\omega_{\bf t})]_{e},
$$
where $\circ$ denotes composition.
\end{lemma}

We write $\mathfrak{A}\stackrel{\sim}{\rightarrow} \mathfrak{B}$ to indicate that there is a bijection $\Omega$ between $\mathfrak{A}$ and $\mathfrak{B}$. If $\omega$ is a weight function of $\mathfrak{A}$, we say that $\Omega$ is $\omega$-preserving if there exists a weight function $\tau$ of $\mathfrak{B}$ such that $\omega = \tau\Omega$. In this case, $[(\mathfrak{A},\omega)]_{e} = [(\mathfrak{B},\tau)]_{e}$. If for sets $\mathfrak{A}, \mathfrak{B}, \mathfrak{C}$ and a mapping $\Omega$ we have $\Omega: \mathfrak{A}\stackrel{\sim}{\rightarrow} \mathfrak{B}\star \mathfrak{C}$, then we say that {i)} $\Omega$ is a direct decomposition for $\mathfrak{A}$ and {ii)} $\Omega$ is an indirect decomposition for $\mathfrak{B}$ or $\mathfrak{C}$ (since any of these sets are embedded on the RHS). This terminology also applies for composition $\circledast$.

We illustrate these ideas with the following brief examples.

\begin{example}\label{ex1}
Let $\mathfrak{D}$ be the set of all permutations with no fixed
points. We give two approaches to determine the number of such
permutations on $n$ points: one by an indirect decomposition and the
other by a direct decomposition. Let $\mathfrak{P}$ be the set of
all permutations and $\mathfrak{I}$ be the set of all the identity
permutations. Here the null permutation is regarded as a permutation
in each set. Then $\mathfrak{P} \stackrel{\sim}{\rightarrow}
\mathfrak{I} \star \mathfrak{D}$, an indirect decomposition for
$\mathfrak{D}$. Thus since we have $\mathfrak{P}
\stackrel{\sim}{\rightarrow} \mathfrak{O}$, by the $\star$-Product
Lemma \ref{ProdLemma} and \eqref{elsets}, $(1-x)^{-1} = \exp x \cdot
[(\mathfrak{D},\omega)]_{e}$, where $\omega(\pi)=n$ if $\pi$ is a
permutation on precisely $n$ symbols. This gives the generating
series for $(\mathfrak{D},\omega)$.
\end{example}

\begin{example} \label{ex2}
Alternatively, any permutation can be decomposed into an unordered
set of disjoint cycles. Let $\mathfrak{C}$ be the set of all
canonical non-null cycles. Then $\mathfrak{P}
\stackrel{\sim}{\rightarrow} \mathfrak{U}\circledast \mathfrak{C}$.
By the Composition Lemma, $(1-x)^{-1} = \exp
[(\mathfrak{C},\omega)]_{e}$ where $\omega$ is the same as above,
giving the generating series for $(\mathfrak{C},\omega)$.
Restricting $\mathfrak{P}$ to $\mathfrak{D}$ we have $\mathfrak{D}
\stackrel{\sim}{\rightarrow} \mathfrak{U}\circledast
(\mathfrak{C}-\mathfrak{C}_1)$, where $\mathfrak{C}_1$ is the set of
all $1$-cycles. This is a direct decomposition for $\mathfrak{D}$
and by the Composition Lemma \ref{CompLemma} and \eqref{elsets}
gives $[(\mathfrak{D},\omega)]_{e} = (1-x)^{-1} \exp (-x)$.
\end{example}

\begin{example}\label{ex3}
A more general instance of Example \ref{ex2} is the following well
known result (implicit in the work of Hurwitz) on graph enumeration.
Let $\mathfrak{G}$ be the set of simple vertex-labelled graphs and
$\mathfrak{G}_c$ be the set of simple connected vertex-labelled
graphs, then we have the following decomposition $\mathfrak{G}
\stackrel{\sim}{\rightarrow} \mathfrak{U} \circledast
\mathfrak{G}_c$, since a vertex-labelled graph is an unordered set
of its connected vertex-labelled graphs. Therefore, if $\omega$ is a
weight function that counts the number of vertices then $
[(\mathfrak{G},\omega)]_{e} = \exp [(\mathfrak{G}_{c},\omega)]_{e}$
by the Composition Lemma \ref{CompLemma} and \eqref{elsets}. We
shall adapt this last result to the graphs in $\mathcal{D}(\Y)$.
\end{example}

% 5 -----------------------------------------------------------------------------------------------

\section[Labelled graphs]{Labelled graphs}\label{S:LG}

\subsection[Labelled uni-trivalent graphs]{Labelled uni-trivalent
graphs}

When we apply $\bbr{\cdot}{\cdot}$ to a power series $D$ in
$\mathcal{D}(\Y)$, multiple copies of the same trivalent graph can
come from one or more graphs in $D$. These are required to be
distinguishable and to do this, we introduce labels and decorations
into the graphs.

We label the set of all isomorphic components of a $\Y$-coloured
uni-trivalent diagram ($\Y\neq \emptyset$) by labels from the same
set. The set of labels for different isomorphism classes of
components are pairwise mutually disjoint. We shall use label sets
$\{1^{1}, 2^{1}, \ldots \}$, $\{1^{2}, 2^{2}, \ldots \}, \ldots$ for
this purpose. These labels are applied to components. Moreover, the
univalent vertices with the same colour in a component are
conveniently distinguished by decorating their univalent edges. To
avoid cluttering the diagrams, we have indicated the decoration by
using edges of different thicknesses. We call such diagrams {\em
component labelled diagrams}. Let $\mathcal{D}_{\ell}(\Y)$ be the
algebra over $\mathbb{Q}$ of formal power series in these graphs as
indeterminants. An example of a component labelled diagram is

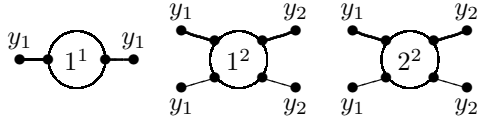
\[
\xy
% shrink
0;/r.30pc/:
%graph 1
(-2,0)*{}="A"; (1,0)*{}="B"; (7,0)*{}="C"; (10,0)*{}="D";
(-2,0)*{\bullet}; (1,0)*{\bullet}; (7,0)*{\bullet};
(10,0)*{\bullet}; "A"; "B" **\dir{-}; "C"; "D"
**\crv{~*=<.7pt>{.}(8,0)}; (4,0)*\xycircle(3,3){-}="f"; (4,0)*{1^{1}};
"f"; (-2,2)*{y_1}; (10,2)*{y_1};
%graph2
(15,5)*{y_1}; (15,-5)*{y_1}; (27,5)*{y_2}; (27,-5)*{y_2};
(15,3)*{}="A"; (18.5,1.9)*{}="B"; (15,-3)*{}="C";
(18.5,-1.9)*{}="D"; (27,3)*{}="E"; (23.5,1.9)*{}="F";
(27,-3)*{}="G"; (23.5,-1.9)*{}="H"; (15,3)*{\bullet};
(18.5,1.9)*{\bullet}; (15,-3)*{\bullet}; (18.5,-1.9)*{\bullet};
(27,3)*{\bullet}; (23.5,1.9)*{\bullet}; (27,-3)*{\bullet};
(23.5,-1.9)*{\bullet}; "A"; "B" **\crv{~*=<0.5pt>{.}(16.75, 2.45)};
"C"; "D" **\dir{-}; "E"; "F" **\crv{~*=<0.5pt>{.}(25.25, 2.45)};
"G"; "H" **\dir{-}; (21,0)*\xycircle(3,3){-}="f"; "f";
(21,0)*{1^{2}};
%\graph3
(33,5)*{y_1}; (33,-5)*{y_1};(45,5)*{y_2}; (45,-5)*{y_2};
(33,3)*{}="A"; (36.5,1.9)*{}="B"; (33,-3)*{}="C";
(36.5,-1.9)*{}="D"; (45,3)*{}="E"; (41.5,1.9)*{}="F";
(45,-3)*{}="G"; (41.5,-1.9)*{}="H"; (33,3)*{\bullet};
(36.5,1.9)*{\bullet}; (33,-3)*{\bullet}; (36.5,-1.9)*{\bullet};
(45,3)*{\bullet}; (41.5,1.9)*{\bullet}; (45,-3)*{\bullet};
(41.5,-1.9)*{\bullet}; "A"; "B" **\crv{~*=<0.5pt>{.}(34.75, 2.45)};
"C"; "D" **\dir{-}; "E"; "F" **\crv{~*=<0.5pt>{.}(43.25, 2.45)};
"G"; "H" **\dir{-}; (39,0)*\xycircle(3,3){-}="f"; "f";
(39,0)*{2^{2}};
\endxy
\]

Note that the components labelled $1^{1}$ and $2^{1}$ are isomorphic
as graphs, so their labels belong to the same label set. In
addition, there are two vertices in the component labelled $1^{1}$
that have the same colour $y_1$, and these are distinguished by
decorating their two incident edges (the two edges have different
thicknesses).

We use $g$ and $h$ to denote  unlabelled graphs and $u_i, v_j$ to
denote unlabelled connected graphs in $\mathcal{D}(\Y)$. To
distinguish between unlabelled and labelled structures, we use $\g$
and $\h$ to denote labelled graphs and ${\sf u}_i$, ${\sf v}_j$ to
denote connected labelled graphs in $\mathcal{D}_{\ell}(\Y)$. The
graph $ {\sf u}_1^{\rho_1} \cdots {\sf u}_n^{\rho_n}$ has labels $\{
1^{1}, 2^{1}, \ldots , \rho_{1}^{1}\}, \ldots , \{1^{n}, 2^{n},
\ldots, \rho_{n}^{n} \}$ for its components ${\sf u}_1, \ldots ,
{\sf u}_n$. To accommodate the labelling we use the exponential
basis of $\mathcal{D}_{\ell}(\Y)$, namely $\{{\sf u}_i^{j}/j ! ~|~
i,j=0,1,2,\ldots \} $. As in $\mathcal{D}(\Y)$, there is only one
graph in $\mathcal{D}_{\ell}(\Y)$ with $\rho_{i}$ components of
${\sf u}_i$ for $1\leq i \leq n$, this is ${\sf u}_1^{\rho_1}\cdots
{\sf u}_n^{\rho_n},$. so there is a one-to-one correspondence
between graphs in $\mathcal{D}_{\ell}(\Y)$ and $\mathcal{D}(\Y)$. In
fact, if $\lambda: \mathfrak{D}(\Y) \rightarrow
\mathfrak{D}_{\ell}(\Y) : u_i^{\rho_i}  \mapsto {\sf u}_i^{\rho_i}$
is the operator that labels elements of $\mathfrak{D}(\Y)$, then it
is an isomorphism. Its inverse $\lambda^{-1}$ discards labels from
the graphs, that is $ \lambda^{-1}:
\mathfrak{D}_{\ell}(\Y)\rightarrow \mathfrak{D}(\Y): {\sf
u}_i^{\rho_i} \mapsto {u}_i^{\rho_i}$.

\subsection[The $\dbbr{\cdot}{\cdot}$ operator]{The $\dbbr{\cdot}{\cdot}$ operator}

We require an operator on $\mathcal{D}_{\ell}(\Y)$ that corresponds
to the operator $\bbr{\cdot}{\cdot}$ in $\mathcal{D}(\Y)$. Such
operator will take graphs $\g, \h$ with labelled components and give
the sum of all the ways of identifying all of the $y$-coloured
vertices in $\g$ with those in $\h$ for all colours $y\in \Y$. The
identification involves the homeomorphic reduction of the path so
formed and the attachment of colour $y$ to the resulting edge. Thus
$\colon \xy (0,0)*{\xy (0,-1.5)*{}="A"; (0,1.5)*{}="B"; (3,1.5)*{}="C";
(3,-1.5)*{}="D";(3,0)*{}="E"; (7,0)*{}="F"; "A"; "B"
**\dir{-};
"A"; "B" **\dir{-}; "B"; "C" **\dir{-}; "C"; "D" **\dir{-}; "D"; "A"
**\dir{-}; "E"; "F" **\dir{-};
"E"*{\bullet}; "F"*{\bullet}; (7,-2)*{y};
\endxy };
(10,0)*{\xy (4,-1.5)*{}="A"; (4,1.5)*{}="B"; (7,1.5)*{}="C";
(7,-1.5)*{}="D";(4,0)*{}="E"; (0,0)*{}="F"; "A"; "B"
**\dir{-};
"A"; "B" **\dir{-}; "B"; "C" **\dir{-}; "C"; "D" **\dir{-}; "D"; "A"
**\dir{-}; "E"; "F" **\dir{-};
"E"*{\bullet}; "F"*{\bullet}; (0,-2)*{y};
\endxy };
(20,0)*{\leadsto}; (32,0)*{\xy (0,-1.5)*{}="A"; (0,1.5)*{}="B";
(3,1.5)*{}="C"; (3,-1.5)*{}="D";  "A"; "B"
**\dir{-}; "B"; "C" **\dir{-}; "C"; "D" **\dir{-}; "D"; "A"
**\dir{-};
(10,-1.5)*{}="A"; (10,1.5)*{}="B"; (13,1.5)*{}="C";
(13,-1.5)*{}="D"; "A"; "B"
**\dir{-}; "B"; "C" **\dir{-}; "C"; "D" **\dir{-}; "D"; "A"
**\dir{-};
(3,0)*{}="E"; (10,0)*{}="F"; "E"; "F" **\dir{-}; "E"*{\bullet};
"F"*{\bullet}; (6.5,-2)*{y};
\endxy };
\endxy
$, where $\xy (0,-1.5)*{}="A"; (0,1.5)*{}="B"; (3,1.5)*{}="C";
(3,-1.5)*{}="D";(3,0)*{}="E"; (7,0)*{}="F"; "A"; "B"
**\dir{-};
"A"; "B" **\dir{-}; "B"; "C" **\dir{-}; "C"; "D" **\dir{-}; "D"; "A"
**\dir{-}; "E"; "F" **\dir{-};
"E"*{\bullet}; "F"*{\bullet}; (7,-2)*{y}; \endxy
$ denotes a
uni-trivalent diagram with a specified $y$-coloured vertex. We keep
decoration in edges when their end vertices are identified
(thus, half of an edge may have a different thickness than the other half) . This
gives trivalent graphs with some coloured edges, some decorated
edges and with the uni-trivalent components as labelled subgraphs.
We denote the algebra of power series in such graphs by
$\mathcal{D}_{\ell}(\emptyset)$. The following is an instance of
such a graph, with the subgraphs highlighted by dotted closed
curves:

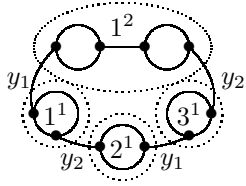
\[
\xy
% shrink
0;/r.35pc/: (2,0)*\xycircle(2,2){-}; (10,0)*\xycircle(2,2){-};
(6,-9)*\xycircle(2,2){-}; (0,-6)*\xycircle(2,2){-};
(12,-6)*\xycircle(2,2){-}; (6,2)*{1^{2}}; (6,-9)*{2^{1}};
(0,-6)*{1^{1}}; (12,-6)*{3^{1}}; (4,0)*{}="X1"; (8,0)*{}="X2";
"X1"*{\bullet}; "X2"*{\bullet}; "X1";"X2" **\dir{-}; (0,0)*{}="X1";
(-2,-6)*{}="X2"; "X1"*{\bullet}; "X2"*{\bullet}; "X1";"X2"
**\crv{(-3,-3)}; (12,0)*{}="X1"; (14,-6)*{}="X2"; "X1"*{\bullet};
"X2"*{\bullet}; "X1";"X2" **\crv{(15,-3)};  (12,-8)*{}="X1";
(8,-9)*{}="X2";
 "X1"*{\bullet};
"X2"*{\bullet}; "X1";"X2" **\crv{(10,-9)}; (0,-8)*{}="X1";
(4,-9)*{}="X2"; "X1"*{\bullet}; "X2"*{\bullet}; "X1";"X2"
**\crv{(2,-9)};  (-3.3,-3)*{y_1}; (16,-3)*{y_2}; (10.5,-10.3)*{y_1};
(1.5,-10.3)*{y_2}; (6,-9)*\xycircle(3,3){.}="g1";
(12,-6)*\xycircle(3,3){.}="g2"; (0,-6)*\xycircle(3,3){.}="g3";
(6,-0)*\xycircle(8,4){.}="g4";
\endxy
\]

We define such operator in the exponential basis of
$\mathcal{D}_{\ell}(\Y)$, and extend it bilinearly and denote it by
$\dbbr{\cdot}{\cdot}:\mathcal{D}_{\ell}(\Y) \otimes
\mathcal{D}_{\ell}(\Y) \rightarrow \mathcal{D}_{\ell}(\emptyset)$.
Thus $\dbbrc{\cdot}{\cdot}$ is the connected part of
$\dbbr{\cdot}{\cdot}$. For the above combinatorial reason we use
label sets $\{1^{1}_{l}, 2^{1}_{l},\ldots\}, \{1^{2}_l,
2^{2}_l,\ldots\},\ldots$ for components in the left argument of
$\dbbr{\cdot}{\cdot}$ and $\{1^{1}_{r}, 2^{1}_{r},\ldots\},
\{1^{2}_r, 2^{2}_r,\ldots\},\ldots$ for the components in the right
argument of $\dbbr{\cdot}{\cdot}$. What we have gained with the
component labelling and decorations is that $\dbbr{\cdot}{\cdot}$
will give a sum of different trivalent diagrams. Moreover, any
rational coefficients in the trivalent diagrams come directly from
the original graphs in $\mathcal{D}_{\ell}(\Y)$.

% presentation graph 1
\[
\xy
% shrink
0;/r.35pc/: (-4,0)*{\Big\langle \Big\langle}; (3,0)*{\xy
(0,0)*{}="A"; (2,0)*{}="B";
          (6,0)*{}="C"; (8,0)*{}="D";
      "A"*{\bullet}; "B"*{\bullet}; "C"*{\bullet}; "D"*{\bullet};
      "A"; "B" **\dir{-}; "C"; "D" **\crv{~*=<0.5pt>{.}(7,0)};
      (4,0)*\xycircle(2,2){-}="f"; (4,0)*{1^{1}_l};
    "f";
    (1,-2)*{y_1};
    (8,-2)*{y_2};
    \endxy
    };
(19,0)*{\xy
    (0,0)*{}="A"; (2,0)*{}="B"; (6,0)*{}="C"; (10,0)*{}="D";
    (14,0)*{}="E"; (16,0)*{}="F"; "A"; "B" **\dir{-}; "C"; "D"
    **\dir{-}; "E"; "F" **\dir{-}; (4,0)*\xycircle(2,2){-}="f1";
    (12,0)*\xycircle(2,2){-}="f2"; (8,-2)*{1^{2}_l}; "f1"; "f2";
    (1,-2)*{y_1};
    (16,-2)*{y_2};
    "A"*{\bullet}; "B"*{\bullet};"C"*{\bullet}; "D"*{\bullet};
     "E"*{\bullet}; "F"*{\bullet};
    \endxy
    };
(30,-1)*{,}; (37,0)*{\xy (0,0)*{}="A";
    (2,0)*{}="B";
    (6,0)*{}="C"; (8,0)*{}="D";
    "A"*{\bullet}; "B"*{\bullet}; "C"*{\bullet}; "D"*{\bullet};
    "A"; "B" **\dir{-}; "C"; "D" **\dir{-};
    (4,0)*\xycircle(2,2){-}="f"; "f"; (4,0)*{1^{1}_r};
    (1,-2)*{y_1};
    (8,-2)*{y_2};
    \endxy
    };
(49,0)*{\xy (0,0)*{}="A"; (2,0)*{}="B";
    (6,0)*{}="C"; (8,0)*{}="D"; "A"; "B" **\dir{-};
     "A"*{\bullet}; "B"*{\bullet}; "C"*{\bullet}; "D"*{\bullet};
    "C"; "D" **\dir{-};
    (4,0)*\xycircle(2,2){-}="f"; "f"; (4,0)*{2^{1}_r};
    (1,-2)*{y_1};
    (8,-2)*{y_2};
    \endxy
    };
(56,0)*{\Big\rangle \Big\rangle}; (60,0)*{=}; (68,0)*{\xy
    (-5,0)*{}="A"; (-1,0)*{}="B";
    (-5,-6)*{}="C"; (-1,-6)*{}="D";
    (-5,0)*{\bullet}; (-1,0)*{\bullet};
    (-5,-6)*{\bullet}; (-1,-6)*{\bullet};
    (-3,0)*\xycircle(2,2){-}="f1";
    (-3,-6)*\xycircle(2,2){-}="f2";
    "A";"C" **\crv{(-7,-3)};
    "B";"D" **\crv{(1,-3)};
    (-7,-3)*{y_1};
    (1.4,-3)*{y_2};
    "f1"; "f2";
    (-3,0)*{1^{1}_{l}};
    (-3,-6)*{1^{1}_{r}};
    \endxy
       };
(81,1)*{ \xy (2,0)*\xycircle(2,2){-}; (10,0)*\xycircle(2,2){-};
(6,-7)*\xycircle(2,2){-}; (6,2)*{1^{2}_l}; (6,-7)*{2^{1}_r};
(4,0)*{}="X1"; (8,0)*{}="X2";
 "X1"*{\bullet}; "X2"*{\bullet};
"X1";"X2" **\dir{-};
%"X1";"X2" **\crv{(6,4)};
(0,0)*{}="X1";
(4,-7)*{}="X2";
 "X1"*{\bullet}; "X2"*{\bullet};
"X1";"X2" **\crv{(1,-3) & (2,-5) };
(12,0)*{}="X1";
(8,-7)*{}="X2";
 "X1"*{\bullet}; "X2"*{\bullet};
"X1";"X2" **\crv{(11,-3) & (10,-5) }; (0.6,-4)*{y_1};
(11.9,-4)*{y_2};
\endxy
};
%equation label
(-13,-9)*{(3)};
%line 2
(-2,-15)*{+}; (6,-15)*{\xy
    (-5,0)*{}="A"; (-1,0)*{}="B";
    (-5,-6)*{}="C"; (-1,-6)*{}="D";
    (-5,0)*{\bullet}; (-1,0)*{\bullet};
    (-5,-6)*{\bullet}; (-1,-6)*{\bullet};
    (-3,0)*\xycircle(2,2){-}="f1";
    (-3,-6)*\xycircle(2,2){-}="f2";
    "A";"C" **\crv{(-7,-3)};
    "B";"D" **\crv{(1,-3)};
    "f1"; "f2";
    (-3,0)*{1^{1}_{l}};
    (-3,-6)*{2^{1}_{r}};
    (-7,-3)*{y_1};
    (1.4,-3)*{y_2};
    \endxy
       };
(19,-14)*{ \xy (2,0)*\xycircle(2,2){-}; (10,0)*\xycircle(2,2){-};
(6,-7)*\xycircle(2,2){-}; (6,2)*{1^{2}_l}; (6,-7)*{1^{1}_r};
(4,0)*{}="X1"; (8,0)*{}="X2";
 "X1"*{\bullet}; "X2"*{\bullet};
"X1";"X2" **\dir{-};
%"X1";"X2" **\crv{(6,4)};
(0,0)*{}="X1"; (4,-7)*{}="X2";
 "X1"*{\bullet}; "X2"*{\bullet};
"X1";"X2" **\crv{(1,-3) & (2,-5) }; (12,0)*{}="X1"; (8,-7)*{}="X2";
 "X1"*{\bullet}; "X2"*{\bullet};
"X1";"X2" **\crv{(11,-3) & (10,-5) }; (0.6,-4)*{y_1};
(11.9,-4)*{y_2};
\endxy
}; (29,-15)*{+}; (43,-15)*{ \xy (2,0)*\xycircle(2,2){-};
(10,0)*\xycircle(2,2){-}; (6,-9)*\xycircle(2,2){-};
(0,-6)*\xycircle(2,2){-}; (12,-6)*\xycircle(2,2){-};
(6,2)*{1^{2}_l}; (6,-9)*{1^{1}_l}; (0,-6)*{2^{1}_r};
(12,-6)*{1^{1}_r}; (4,0)*{}="X1"; (8,0)*{}="X2";
 "X1"*{\bullet}; "X2"*{\bullet};
"X1";"X2" **\dir{-}; (0,0)*{}="X1"; (-2,-6)*{}="X2";
 "X1"*{\bullet}; "X2"*{\bullet};
"X1";"X2" **\crv{(-3,-3)}; (12,0)*{}="X1"; (14,-6)*{}="X2";
 "X1"*{\bullet}; "X2"*{\bullet};
"X1";"X2" **\crv{(15,-3)}; (12,-8)*{}="X1"; (8,-9)*{}="X2";
 "X1"*{\bullet}; "X2"*{\bullet};
"X1";"X2" **\crv{(10,-9)}; (0,-8)*{}="X1"; (4,-9)*{}="X2";
"X1"*{\bullet}; "X2"*{\bullet}; "X1";"X2"
**\crv{(2,-9)};
(-3,-3)*{y_1}; (16,-3)*{y_2}; (11,-10)*{y_1}; (1,-10)*{y_2};
\endxy
}; (56,-15)*{+}; (70,-15)*{ \xy (2,0)*\xycircle(2,2){-};
(10,0)*\xycircle(2,2){-}; (6,-9)*\xycircle(2,2){-};
(0,-6)*\xycircle(2,2){-}; (12,-6)*\xycircle(2,2){-};
(6,2)*{1^{2}_l}; (6,-9)*{1^{1}_l}; (0,-6)*{1^{1}_r};
(12,-6)*{2^{1}_r}; (4,0)*{}="X1"; (8,0)*{}="X2";
 "X1"*{\bullet}; "X2"*{\bullet};
"X1";"X2" **\dir{-};
(0,0)*{}="X1";
(-2,-6)*{}="X2";
 "X1"*{\bullet}; "X2"*{\bullet};
"X1";"X2" **\crv{(-3,-3)};
(12,0)*{}="X1";
(14,-6)*{}="X2";
 "X1"*{\bullet}; "X2"*{\bullet};
"X1";"X2" **\crv{(15,-3)};
(12,-8)*{}="X1";
(8,-9)*{}="X2";
 "X1"*{\bullet}; "X2"*{\bullet};
"X1";"X2" **\crv{(10,-9)};
(0,-8)*{}="X1";
(4,-9)*{}="X2";
 "X1"*{\bullet}; "X2"*{\bullet};
"X1";"X2" **\crv{(2,-9)}; (-3,-3)*{y_1}; (16,-3)*{y_2};
(11,-10)*{y_1}; (1,-10)*{y_2};
\endxy
};
\endxy
\]

Note that the first two graphs in the right hand side have different
subgraph labelling and the last two are distinguishable by their
edge colouring.

As in the case of $\mathcal{D}(\Y)$, where $\Y \neq \emptyset$, let
$G$ and $\Gamma_{i}$ denote a trivalent graph and a trivalent
connected graph in $\mathcal{D}(\emptyset)$. Let ${\sf G}$ and ${\sf
\Gamma}_i$ denote the labelled and decorated counterparts in
$\mathcal{D}_{\ell}(\emptyset)$. Again, if $\lambda^{-1}:
\mathfrak{D}_{\ell}(\emptyset)\rightarrow \mathfrak{D}(\emptyset)$
discards subgraph labels, edge colourings and decorations, then
$\lambda^{-1}: {\sf \Gamma}_i \mapsto \Gamma_i$. In this case
$\lambda^{-1}$ is not necessarily one-to-one.

% 6 -----------------------------------------------------------------------------------------------

\section[The graph-subgraph exponential generating series]{The graph-subgraph exponential generating series}\label{S:GSEGS}

\subsection[Component-Subgraph Decomposition]{Component-Subgraph
Decomposition}

To account for the appearance of the external exponential function
in the right hand side of Theorem \ref{thm2}, we show here that a
graph in $\mathcal{D}_{\ell}(\emptyset)$ can be decomposed into an
unordered set of its components and that this decomposition
preserves the subgraph labelling. The occurrences of the remaining
two internal exponential functions in the right hand side of Theorem
\ref{thm2} will be explained very simply in Section 7.

For simplicity, let ${\sf G}= \prod_{k=1}^{m} {\sf \Gamma}_{k} \in
\mathcal{D}_{\ell}(\emptyset)$ be one of the terms of  $\dbbr{ {\sf
u}^{\rho}}{{\sf v}^{\varrho}}$. ${\sf G}$ has the subgraphs ${\sf
u}$ with all the labels $\{1_{l}^{1}, 2_{l}^{1},\ldots,
\rho_{l}^{1}\}$. Each of the connected components ${\sf \Gamma_k}$
has $\rho_k$ of these with the labels $\{\alpha^{1}_{(k,1)},
\alpha^{1}_{(k,2)}, \ldots, \alpha^{1}_{(k,\rho_k)}\}$.  Let
$\alpha_{(k)} = \{ \alpha_{(k,1)},\alpha_{(k,2)}, \ldots,
\alpha_{(k,\rho_k)}\}$, without loss of generality $\alpha_{(k,1)} <
\ldots <\alpha_{(k,\rho_k)}$. Similarly, for the subgraph ${\sf v}$
with labels $\{1_{r}^{1}, 2_{r}^{1},\ldots , \varrho_{r}^{1} \}$,
for each component ${\sf \Gamma_k}$, we get a set $\beta_{(k)}=\{
\beta^{1}_{(k,1)}, \beta^{1}_{(k,2)},\ldots,
\beta^{1}_{(k,\varrho_k)}\}$.

Since ${\sf \Gamma}_k \in \mathcal{D}_{\ell}(\emptyset)$ is the
corresponding graph by replacing the label sets $\alpha_{(k)}$ and
$\beta_{(k)}$ by $\{1,\ldots, \rho_k\}$ and $\{1,\ldots,
\varrho_k\},$ respectively, we can denote the ${\sf
\Gamma}_k$-component of ${\sf G}$ using the notation of Section
\ref{lemmas} as $ _{\alpha_{(k)}}({\sf \Gamma}_k)_{\beta_{(k)}}$.
Hence,
$$
\prod_{k=1}^{m} {\sf \Gamma}_{k} \longmapsto \mathsf{X}_{k=1}^{m}~_{
\alpha_{(k)}}({\sf \Gamma}_k)_{\beta_{(k)}} \text{ \quad where
\quad}  {\biguplus}_{k=1}^{m} \alpha_{(k)} = \mathcal{N}_{\rho}
\text{~and ~}  {\biguplus}_{k=1}^{m} \beta_{(k)} =
\mathcal{N}_{\varrho}.
$$
Where $\biguplus$ indicates that the sets of labels $\alpha_{(k)}$ for $k=1,\ldots, m$ are mutually
disjoint. This also holds for the sets $\beta_{(k)}$. For example,
from $(3)$,

% graph 6
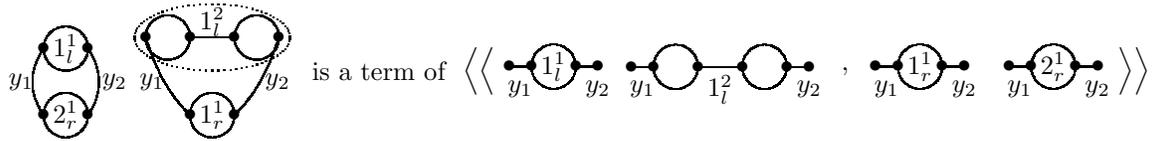
\[
\xy 0;/r.35pc/: (0,-1)*{\xy
    (-5,0)*{}="A"; (-1,0)*{}="B";
    (-5,-6)*{}="C"; (-1,-6)*{}="D";
    (-5,0)*{\bullet}; (-1,0)*{\bullet};
    (-5,-6)*{\bullet}; (-1,-6)*{\bullet};
    (-3,0)*\xycircle(2,2){-}="f1";
    (-3,-6)*\xycircle(2,2){-}="f2";
    "A";"C" **\crv{(-7,-3)};
    "B";"D" **\crv{(1,-3)};
    "f1"; "f2";
    (-3,0)*{1^{1}_{l}};
    (-3,-6)*{2^{1}_{r}};
    (-7,-3)*{y_1};
    (1.4,-3)*{y_2};
    \endxy
       };
%next graph
(13,0)*{ \xy (2,0)*\xycircle(2,2){-}; (10,0)*\xycircle(2,2){-};
(6,-7)*\xycircle(2,2){-}; (6,0)*\xycircle(7,3){.};
(6,1.4)*{1^{2}_{l}}; (6,-7)*{1^{1}_{r}}; (4,0)*{}="X1";
(8,0)*{}="X2";
 "X1"*{\bullet}; "X2"*{\bullet};
"X1";"X2" **\dir{-};
%"X1";"X2" **\crv{(6,4)};
(0,0)*{}="X1";
(4,-7)*{}="X2";
 "X1"*{\bullet}; "X2"*{\bullet};
"X1";"X2" **\crv{(1,-3) & (2,-5) };
(12,0)*{}="X1";
(8,-7)*{}="X2";
 "X1"*{\bullet}; "X2"*{\bullet};
"X1";"X2" **\crv{(11,-3) & (10,-5) }; (0.6,-4)*{y_1};
(11.9,-4)*{y_2};
\endxy
};
%sign
(28,0)*{\text{is a term of}};
%nextgraph
(37,0)*{\Big\langle \Big\langle}; (44,0)*{\xy (0,0)*{}="A";
(2,0)*{}="B";
          (6,0)*{}="C"; (8,0)*{}="D";
      "A"*{\bullet}; "B"*{\bullet}; "C"*{\bullet}; "D"*{\bullet};
      "A"; "B" **\dir{-}; "C"; "D" **\dir{-};
      (4,0)*\xycircle(2,2){-}="f"; (4,0)*{1^{1}_{l}};
    "f";
    (1,-2)*{y_1};
    (8,-2)*{y_2};
    \endxy
    };
(59,-0.3)*{\xy
    (0,0)*{}="A"; (2,0)*{}="B"; (6,0)*{}="C"; (10,0)*{}="D";
    (14,0)*{}="E"; (16,0)*{}="F"; "A"; "B" **\dir{-}; "C"; "D"
    **\dir{-}; "E"; "F" **\dir{-}; (4,0)*\xycircle(2,2){-}="f1";
    (12,0)*\xycircle(2,2){-}="f2"; (8,-2)*{1^{2}_{l}}; "f1"; "f2";
    (1,-2)*{y_1};
    (16,-2)*{y_2};
    "A"*{\bullet}; "B"*{\bullet};"C"*{\bullet}; "D"*{\bullet};
     "E"*{\bullet}; "F"*{\bullet};
    \endxy
    };
(70,0)*{{\bf ,}}; (77,0)*{\xy (0,0)*{}="A";
    (2,0)*{}="B";
    (6,0)*{}="C"; (8,0)*{}="D";
    "A"*{\bullet}; "B"*{\bullet}; "C"*{\bullet}; "D"*{\bullet};
    "A"; "B" **\dir{-}; "C"; "D" **\dir{-};
    (4,0)*\xycircle(2,2){-}="f"; "f"; (4,0)*{1^{1}_{r}};
    (1,-2)*{y_1};
    (8,-2)*{y_2};
    \endxy
    };
(89,0)*{\xy (0,0)*{}="A"; (2,0)*{}="B";
    (6,0)*{}="C"; (8,0)*{}="D"; "A"; "B" **\dir{-};
     "A"*{\bullet}; "B"*{\bullet}; "C"*{\bullet}; "D"*{\bullet};
    "C"; "D" **\dir{-};
    (4,0)*\xycircle(2,2){-}="f"; "f"; (4,0)*{2^{1}_{r}};
    (1,-2)*{y_1};
    (8,-2)*{y_2};
    \endxy
    };
(96,0)*{\Big\rangle \Big\rangle};
\endxy
\]

\noindent and it can be decomposed into

% graph 7
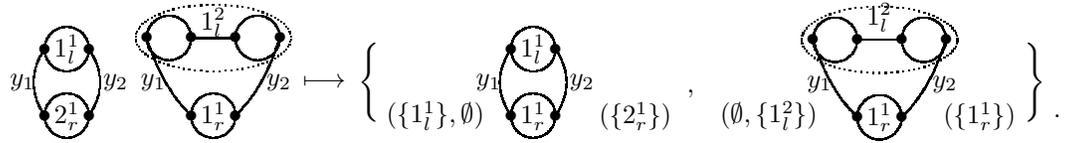
\[
\xy 0;/r.35pc/:  (0,0)*{\xy
    (-5,0)*{}="A"; (-1,0)*{}="B";
    (-5,-6)*{}="C"; (-1,-6)*{}="D";
    (-5,0)*{\bullet}; (-1,0)*{\bullet};
    (-5,-6)*{\bullet}; (-1,-6)*{\bullet};
    (-3,0)*\xycircle(2,2){-}="f1";
    (-3,-6)*\xycircle(2,2){-}="f2";
    "A";"C" **\crv{(-7,-3)};
    "B";"D" **\crv{(1,-3)};
    "f1"; "f2";
    (-3,0)*{1^{1}_{l}};
    (-3,-6)*{2^{1}_{r}};
    (-7,-3)*{y_1};
    (1.4,-3)*{y_2};
    \endxy
       };
%next graph
(13,1)*{ \xy (2,0)*\xycircle(2,2){-}; (10,0)*\xycircle(2,2){-};
(6,-7)*\xycircle(2,2){-}; (6,0)*\xycircle(7,3){.};
(6,1.4)*{1^{2}_{l}}; (6,-7)*{1^{1}_{r}}; (4,0)*{}="X1";
(8,0)*{}="X2";
 "X1"*{\bullet}; "X2"*{\bullet};
"X1";"X2" **\dir{-};
%"X1";"X2" **\crv{(6,4)};
(0,0)*{}="X1";
(4,-7)*{}="X2";
 "X1"*{\bullet}; "X2"*{\bullet};
"X1";"X2" **\crv{(1,-3) & (2,-5) }; (12,0)*{}="X1"; (8,-7)*{}="X2";
 "X1"*{\bullet}; "X2"*{\bullet};
"X1";"X2" **\crv{(11,-3) & (10,-5) }; (0.6,-4)*{y_1};
(11.9,-4)*{y_2};
\endxy
};
%sing
(23,0)*{\longmapsto};
%nextgraph
(27,0)*{\Biggl\{}; (33,-3)*{(\{1^{1}_l\}, \emptyset)}; (42,0)*{\xy
    (-5,0)*{}="A"; (-1,0)*{}="B";
    (-5,-6)*{}="C"; (-1,-6)*{}="D";
    (-5,0)*{\bullet}; (-1,0)*{\bullet};
    (-5,-6)*{\bullet}; (-1,-6)*{\bullet};
    (-3,0)*\xycircle(2,2){-}="f1";
    (-3,-6)*\xycircle(2,2){-}="f2";
    "A";"C" **\crv{(-7,-3)};
    "B";"D" **\crv{(1,-3)};
    "f1"; "f2";
    (-3,0)*{1^{1}_{l}};
    (-3,-6)*{1^{1}_{r}};
    (-7,-3)*{y_1};
    (1.4,-3)*{y_2};
    \endxy
       };
(51,-3)*{(\{2^{1}_r\})};
%nextgraph
(56,-1)*{,}; (63,-3)*{(\emptyset,\{1^{2}_l\})}; (73,1)*{ \xy
(2,0)*\xycircle(2,2){-}; (10,0)*\xycircle(2,2){-};
(6,-7)*\xycircle(2,2){-}; (6,0)*\xycircle(7,3){.};
(6,2)*{1^{2}_{l}}; (6,-7)*{1^{1}_{r}}; (4,0)*{}="X1"; (8,0)*{}="X2";
"X1"*{\bullet}; "X2"*{\bullet}; "X1";"X2" **\dir{-};
%"X1";"X2" **\crv{(6,4)};
(0,0)*{}="X1"; (4,-7)*{}="X2"; "X1"*{\bullet}; "X2"*{\bullet};
"X1";"X2" **\crv{(1,-3) & (2,-5) }; (12,0)*{}="X1"; (8,-7)*{}="X2";
"X1"*{\bullet}; "X2"*{\bullet}; "X1";"X2" **\crv{(11,-3) & (10,-5)
}; (0.6,-4)*{y_1}; (11.9,-4)*{y_2};
\endxy
}; (82,-3)*{(\{1^{1}_r \})}; (87,0)*{\Biggr\}}; (89,-3)*{.};
\endxy
\]

\noindent where $\{1_{l}^{1}\} ~{\uplus}~ \emptyset =
\mathcal{N}_1$, $\emptyset ~{\uplus}~ \{1_{l}^{2}\} =
\mathcal{N}_1$ and $\{2_{r}^{1} \} ~{\uplus}~ \{1_{r}^{1} \}=
\mathcal{N}_2$. As noted in the previous section, the rational
coefficients of the trivalent graphs in $\dbbr{\cdot}{\cdot}$ come
directly from the graphs in $\mathcal{D}_{\ell}(\Y)$, so we can
weight the graph with these rational coefficients. This suggests the
following lemma.

\begin{lemma}\label{graphs}
For ${\sf D}_1, {\sf D}_2 \in \mathcal{D}_{\ell}(\Y)$, let
$\mathfrak{G}_{c}({\sf D}_1,{\sf D}_2)$ be the set of
subgraph-labelled connected graphs in $\dbbr{{\sf D}_1}{{\sf D}_2}$
with rational weights and let $\mathfrak{G}({\sf D}_1, {\sf D}_2)$
the set of subgraph-labelled graphs with components from
$\mathfrak{G}_{c}({\sf D}_1,{\sf D}_2)$ (their rational weights are
obtained by multiplying the weights of the components). Then
$$
\mathfrak{G}({\sf D}_1,{\sf D}_2) \stackrel{\sim}{\longrightarrow}
\mathfrak{U} \circledast \mathfrak{G}_{c}({\sf D}_1,{\sf D}_2).
$$
\end{lemma}

\begin{proof}
This is simply an enrichment of Example \ref{ex3}.
\end{proof}

\subsection[The weight functions for the graphs]{The weight function for the graphs}\label{weights}

Let ${\sf G}=\prod_{k=1}^{m} {\sf \Gamma}_k \in
\mathcal{D}_{\ell}(\emptyset)$. Now, for $\Gamma \in
\mathcal{D}(\emptyset)$ and $u \in \mathcal{D}(\Y)$, we define the
following weight functions:

\begin{itemize}
\item[(i)] {\em trivalent component weights:} Let $\omega_{\Gamma}:\mathcal{D}_{\ell}(\emptyset)\rightarrow \mathcal{N}_{\geq 0}$, where $\omega_{\Gamma}(\sf G)$
is the number of components of $\lambda^{-1}({\sf
G})=\prod^{m}_{k=1} \Gamma_k$ that are isomorphic to $\Gamma$.
\item[(ii)] {\em subgraph weights:} Let $\omega_{u}:\mathcal{D}_{\ell}(\emptyset)\rightarrow \mathcal{N}_{\geq 0}$, where $\omega_{u}(\sf
G)$ is the number of appearances of $\lambda(u)={\sf u}$ as a
subgraph in ${\sf G}$.
\item[(iii)] Let $\omega =   \left(\bigotimes_{\Gamma \in
\mathcal{D}(\emptyset)} \omega_{\Gamma} \right) \otimes
\left(\bigotimes_{u \in \mathcal{D}(\Y)} \omega_{u} \right)$.
\end{itemize}
\noindent Similarly, if $\g= \prod_{i=1}^{n} {\sf u}_i \in
\mathcal{D}_{\ell}(\Y)$, we define
\begin{itemize}
\item[(iv)] {\em uni-trivalent component weights:} Let $\theta_u: \mathcal{D}_{\ell}(\Y) \rightarrow \mathcal{N}_{\geq 0}$,
where $\theta_{u}(\sf g)$ is the number of components of
$\lambda^{-1}({\sf g})= \prod_{i=1}^{n} u_i$ that are isomorphic to
$u$.
\item[(v)] Let $\theta = \bigotimes_{u \in \mathcal{D}(\Y)}
\theta_u$.
\end{itemize}

\subsection[The linear operators $\Phi$ and $\Phi^{\circ}$]{The linear functions $\Phi$ and $\Phi^{\circ}$}

We use the linear operators $\Phi$ and $\Phi^{\circ}$, respectively,
for computing the generating series for labelled trivalent graphs
and labelled uni-trivalent graphs. Using $\Gamma \in
\mathcal{D}(\emptyset)$, and $u\in \mathcal{D}(\Y)$ as
indeterminates, $\Gamma$ marks ordinarily the
 number of $\Gamma$-components, and  $u$ marks exponentially the number of $u$-subgraphs. We therefore work
 in the ring of formal power series $ \mathcal{R} \equiv \mathbb{Q}[\prod_{\Gamma \in
\mathcal{D}(\emptyset)}\Gamma ]\,[[\prod_{u \in \mathcal{D}(\Y)} u
]]$ and its subring $\mathcal{R}_{\mathcal{D}(\Y)} \equiv
\mathbb{Q}[[\prod_{u \in \mathcal{D}(\Y)} u ]]$. Let $\omega$ be as
defined in Section \ref{weights}. Then
$$
\Phi:\mathcal{D}_{\ell}(\emptyset)  \rightarrow \mathcal{R}:{\sf G}
\mapsto \left[({\sf G}, \omega)\right]_{(o;e)},
$$
where $\Phi$ is extended linearly to
$\mathcal{D}_{\ell}(\emptyset)$.

As an example, if $\prod_{k=1}^{m} {\sf \Gamma}_{k}$ is a term in $
\dbbr{ {\sf u}_1^{\rho_1}\cdots {\sf u}_{n_1}^{\rho_{n_1}}}{{\sf
v}_1^{\varrho_1}\cdots {\sf v}_{n_2}^{\varrho_{n_2}}}$, then
$$
\Phi(\prod_{k=1}^{m} {\sf \Gamma}_k ) = [(\prod_{k=1}^{m} {\sf
\Gamma}_k,\omega )]_{(o;e)}\!(\Gamma,u) = \left( \prod_{k=1}^{m}
\Gamma_k\right)\left( \prod_{j=1}^{n_1}\frac{u_j^{\rho_j}}{\rho_j!}
\prod_{j=1}^{n_2} \frac{v_j^{\varrho_j}}{\varrho_j!} \right).
$$
Thus, the generating series encodes the trivalent graph without
labels and decorations but with subgraph information. For example,
for the first term in the right hand side of $(3)$, we have
\[
\xy 0;/r.35pc/: (-7,0)*{\Phi :}; (1,-1)*{\xy
    (-5,0)*{}="A"; (-1,0)*{}="B";
    (-5,-6)*{}="C"; (-1,-6)*{}="D";
    (-5,0)*{\bullet}; (-1,0)*{\bullet};
    (-5,-6)*{\bullet}; (-1,-6)*{\bullet};
    (-3,0)*\xycircle(2,2){-}="f1";
    (-3,-6)*\xycircle(2,2){-}="f2";
    "A";"C" **\crv{(-7,-3)};
    "B";"D" **\crv{(1,-3)};
    "f1"; "f2";
    (-3,0)*{1^{1}_{l}};
    (-3,-6)*{2^{1}_{r}};
    (-7,-3)*{y_1};
    (1.4,-3)*{y_2};
    \endxy
       };
(14,0)*{ \xy (2,0)*\xycircle(2,2){-}; (10,0)*\xycircle(2,2){-};
(6,-7)*\xycircle(2,2){-}; (6,2)*{1^{2}_{l}}; (6,-7)*{2^{1}_{r}};
(4,0)*{}="X1"; (8,0)*{}="X2";  (4,0)*{\bullet};
(8,0)*{\bullet};"X1";"X2" **\dir{-};
%"X1";"X2" **\crv{(6,4)};
(0,0)*{}="X1"; (4,-7)*{}="X2"; (0,0)*{\bullet}; (4,-7)*{\bullet};
"X1";"X2" **\crv{(1,-3) & (2,-5) }; (12,0)*{}="X1"; (8,-7)*{}="X2";
"X1";"X2" **\crv{(11,-3) & (10,-5) }; (12,0)*{\bullet};
(8,-7)*{\bullet}; (0.6,-4)*{y_1}; (11.9,-4)*{y_2};
\endxy
};
%sing
(26,0)*{\longmapsto};
%nextgraph
(35,0)*{\xy
    (-5,0)*{}="A"; (-1,0)*{}="B";
    (-5,-6)*{}="C"; (-1,-6)*{}="D";
    (-5,0)*{\bullet}; (-1,0)*{\bullet};
    (-5,-6)*{\bullet}; (-1,-6)*{\bullet};
    (-3,0)*\xycircle(2,2){-}="f1";
    (-3,-6)*\xycircle(2,2){-}="f2";
    "A";"C" **\crv{(-7,-3)};
    "B";"D" **\crv{(1,-3)};
    "f1"; "f2";
     \endxy
       };
(47,0)*{ \xy (2,0)*\xycircle(2,2){-}; (10,0)*\xycircle(2,2){-};
(6,-7)*\xycircle(2,2){-}; (4,0)*{}="X1"; (8,0)*{}="X2";
(4,0)*{\bullet}; (8,0)*{\bullet}; "X1";"X2" **\dir{-};
%"X1";"X2" **\crv{(6,4)};
(0,0)*{}="X1"; (4,-7)*{}="X2"; (0,0)*{\bullet}; (4,-7)*{\bullet};
"X1";"X2" **\crv{(1,-3) & (2,-5) }; (12,0)*{}="X1"; (8,-7)*{}="X2";
"X1";"X2" **\crv{(11,-3) & (10,-5) }; (12,0)*{\bullet};
(8,-7)*{\bullet};
\endxy
}; (60,1)*{\xy (-1,0)*{\Bigl(}; (0,0)*{}="A"; (2,0)*{}="B";
          (6,0)*{}="C"; (8,0)*{}="D"; (9.5,0)*{\Bigr)}; "A"; "B" **\dir{-}; "C"; "D" **\dir{-};
          (0,0)*{\bullet}; (2,0)*{\bullet};
          (6,0)*{\bullet}; (8,0)*{\bullet};
      (4,0)*\xycircle(2,2){-}="f";
    "f";
    (1,-2)*{y_1};
    (8,-2)*{y_2};
    (10.1,2)*{^2};
    (1,-4)*{}="A"; (7,-4)*{}="B";
    "A"; "B" **\dir{-};
    (4,-6)*{2!};
    \endxy
    };
(76,1)*{\xy
    (-1,0)*{\Bigl(}; (0,0)*{}="A"; (2,0)*{}="B"; (6,0)*{}="C";
    (10,0)*{}="D";
    (14,0)*{}="E"; (16,0)*{}="F"; (17.3,0)*{\Bigr)};
    (0,0)*{\bullet}; (2,0)*{\bullet}; (6,0)*{\bullet}; (10,0)*{\bullet};
    (14,0)*{\bullet}; (16,0)*{\bullet};
     "A"; "B" **\dir{-}; "C"; "D"
    **\dir{-}; "E"; "F" **\dir{-}; (4,0)*\xycircle(2,2){-}="f1";
    (12,0)*\xycircle(2,2){-}="f2"; "f1"; "f2";
    (1,-2)*{y_1};
    (16,-2)*{y_2};
    (17.8,2)*{^2};
    (1,-4)*{}="A"; (15,-4)*{}="B";
    "A"; "B" **\dir{-};
    (8,-6)*{2!};
    \endxy
    };
(86,-1)*{.};
\endxy
\]

We also require combinatorial information from elements in
$\mathcal{D}_{\ell}(\Y)$ in the arguments of $\dbbr{\cdot}{\cdot}$.
Let $\theta$ be as defined in Section \ref{weights}. Then
$$
\Phi^{\circ}:\mathcal{D}_{\ell}(\Y)  \rightarrow
\mathcal{R}_{\mathcal{D}(\Y)}: \g \mapsto \left[(\g,
\theta)\right]_{(e)},
$$
where $\Phi^{\circ}$ is extended linearly to
$\mathcal{D}_{\ell}(\Y)$. As an example,
$\Phi^{\circ}(\prod_{j=1}^{n} {\sf u}_j^{\rho_j}) = \prod_{j=1}^{n}
{\sf u}_{j}^{\rho_{j}}/{\rho_{j}!}$ or for the graph in the right
argument of $\dbbr{\cdot}{\cdot}$ in (3),
\[
\xy 0;/r.35pc/: (-10,0)*{\Phi^{\circ} :}; (0,0)*{\xy (0,0)*{}="A";
    (2,0)*{}="B";
    (6,0)*{}="C"; (8,0)*{}="D";
    "A"*{\bullet}; "B"*{\bullet}; "C"*{\bullet}; "D"*{\bullet};
    "A"; "B" **\dir{-}; "C"; "D" **\dir{-};
    (4,0)*\xycircle(2,2){-}="f"; "f"; (4,0)*{1^{1}_{r}};
    (1,-2)*{y_1};
    (8,-2)*{y_2};
    \endxy
    };
(12,0)*{\xy (0,0)*{}="A"; (2,0)*{}="B";
    (6,0)*{}="C"; (8,0)*{}="D"; "A"; "B" **\dir{-};
     "A"*{\bullet}; "B"*{\bullet}; "C"*{\bullet}; "D"*{\bullet};
    "C"; "D" **\dir{-};
    (4,0)*\xycircle(2,2){-}="f"; "f"; (4,0)*{2^{1}_{r}};
    (1,-2)*{y_1};
    (8,-2)*{y_2};
    \endxy
    };
%sing
(23,0)*{\longmapsto};
%nextgraph
(33,1)*{\xy (-1,0)*{\Bigl(}; (0,0)*{}="A"; (2,0)*{}="B";
          (6,0)*{}="C"; (8,0)*{}="D"; (9.5,0)*{\Bigr)}; "A"; "B" **\dir{-}; "C"; "D" **\dir{-};
          (0,0)*{\bullet}; (2,0)*{\bullet};
          (6,0)*{\bullet}; (8,0)*{\bullet};
      (4,0)*\xycircle(2,2){-}="f";
    "f";
    (1,-2)*{y_1};
    (8,-2)*{y_2};
    (10.2,2)*{^2};
    (1,-4)*{}="A"; (7,-4)*{}="B";
    "A"; "B" **\dir{-};
    (4,-6)*{2!};
    \endxy
    };
(40,-1)*{.};
\endxy
\]

If we let $L_{a}^{x} f(x) = f(x)\bigl|_{x=a}$ be the evaluation
operator, the following proposition follows.

\begin{proposition}\label{prop1}
Let $\Y=\{y_1, y_2, \ldots \}$, for ${\sf D}_1, {\sf D}_2 \in \mathcal{D}_{\ell}(\Y)$.
Then
\begin{eqnarray*}
L_{1}^{\mathcal{D}(\Y)} \Phi(\dbbr{{\sf D}_1}{{\sf D}_2}) &=& \bbr{\Phi^{\circ}({\sf D}_1)}{\Phi^{\circ}({\sf D}_2)},\\
L_{1}^{\mathcal{D}(\Y)} \Phi(\dbbrc{{\sf D}_1}{{\sf D}_2}) &=&
\bbrc{\Phi^{\circ}({\sf D}_1)}{\Phi^{\circ}({\sf D}_2)}.
\end{eqnarray*}
\end{proposition}

% 7 -----------------------------------------------------------------------------------------------

\section[Proof of main theorem]{Proof of main theorem}\label{S:PMT}

Before proving of the main theorem (Theorem \ref{thm2}), a
prefactory result is first needed.

\subsection[A prefactory result]{A prefactory result}

Let $B= \sum_{i=0}^{\infty} w_i$ and $C= \sum_{i=0}^{\infty} z_i$
and where $w_i = r_i u_i$ and $z_i = s_i v_i$ are the scalar
products of rational numbers $r_i, s_i$ and connected uni-trivalent
graphs $u_i, v_i \in \mathcal{D}(\Y)$. Let ${\sf B}$ and ${\sf C}$
be the images of $B$ and $C$ in $\mathcal{D}_{\ell}(\Y)$.

Let ${\sf D}_{\sf B}$ be the series of all graphs with components in
${\sf B}$. By a trivial decomposition into components,
$\Phi^{\circ}({\sf D}_{\sf B})= \exp \Phi^{\circ}({\sf B})$. Since
${\sf B}$ is connected then $\Phi^{\circ}({\sf B})= B$ and it
follows that $\Phi^{\circ}({\sf D}_{\sf B}) = \exp B$. Similarly, if
${\sf D}_{\sf C}$ is the series of all graphs with components in
${\sf C}$, then $\Phi^{\circ}({\sf D}_{\sf C}) = \exp C$. Let
$\mathfrak{G}_{c}({\sf D}_{\sf B},{\sf D}_{\sf C})$ and
$\mathfrak{G}({\sf D}_{\sf B},{\sf D}_{\sf C})$ be as in Lemma
\ref{graphs}. We have the following proposition.

\begin{proposition} \label{prop2}
Let ${\sf B}$ and ${\sf C}$ be series of connected graphs in $\mathcal{D}_{\ell}(\Y)$. If ${\sf D}_{\sf B}$ and ${\sf D}_{\sf C}$ are the series of graphs with components from ${\sf B}$ and ${\sf C}$ respectively, then
\begin{equation*}
\label{eqp2} [(\mathfrak{G}({\sf D}_{\sf B}, {\sf D}_{\sf
C}),\omega)]_{(o;e)} = \Phi \dbbr{{\sf D}_{\sf B}}{{\sf D}_{\sf C}}.
\end{equation*}
\end{proposition}

\begin{proof}
We shall show that the terms in $\dbbr{{\sf D}_{\sf B}}{{\sf D}_{\sf
C}}$ are exactly the graphs of $\mathfrak{G}({\sf D}_{\sf B}, {\sf
D}_{\sf C})$. Let ${\sf G}_1$ and ${\sf G}_2$ be terms of
$\dbbr{{\sf D}_{\sf B}}{{\sf D}_{\sf C}}$. More precisely, they come
from $\dbbr{\g_1}{\h_1}$ and $\dbbr{\g_2}{\h_2}$ where $\g_1$ and
$\g_2$ are graphs with components in ${\sf B}$, and $\h_1$ and
$\h_2$ are graphs with components in ${\sf C}$. For any relabelling
of their subgraphs, ${\sf G}_1{\sf G}_2$ is a term of
$\dbbr{\g_1\g_2}{\h_1\h_2}$, where $\g_1\g_2$ and $\h_1\h_2$ also
have components from ${\sf B}$ and ${\sf C}$. Thus ${\sf G}_1{\sf
G}_2$ is a term of $\dbbr{{\sf D}_{\sf B}}{{\sf D}_{\sf C}}$. It
then follows that all the graphs of $\mathfrak{G}({\sf D}_{\sf B},
{\sf D}_{\sf C})$ are terms of $\dbbr{{\sf D}_{\sf B}}{{\sf D}_{\sf
C}}$.

On the other hand, let ${\sf G}_1{\sf G}_2$ be a term of $\dbbr{{\sf
D}_{\sf B}}{{\sf D}_{\sf C}}$. If it is a term in $\dbbr{\g}{\h}$,
we can express $\g = \g_1\g_2$ and $\h=\h_1\h_2$ such that ${\sf
G}_1$ and ${\sf G}_2$ each come from $\dbbr{\g_1}{\h_1}$ and
$\dbbr{\g_2}{\h_2}$. Since  $\g_1$ and $\g_2$ are graphs with
components in ${\sf B}$, and $\h_1$ and $\h_2$ are graphs with
components in ${\sf C}$, ${\sf G}_1$ and ${\sf G}_2$ are terms of
$\dbbr{{\sf D}_{\sf B}}{{\sf D}_{\sf C}}$. This implies that the
connected components of the terms of $\dbbr{{\sf D}_{\sf B}}{{\sf
D}_{\sf C}}$ are all in $\mathfrak{G}_{c}({\sf D}_{\sf B}, {\sf
D}_{\sf C})$. Thus, the terms of $\dbbr{{\sf D}_{\sf B}}{{\sf D}_{\sf
C}}$ are in $\mathfrak{G}({\sf D}_{\sf B},{\sf D}_{\sf C})$. This
gives the desired result, since the weight function $\omega$ is the
same as the one in the definition of $\Phi$.
\end{proof}

\subsection[Proof Theorem \ref{thm2}]{Proof Theorem \ref{thm2}}
We are now in a position to prove the main theorem.

\begin{proof}
Let $B= \sum_{i=0}^{\infty} w_i$ and $C= \sum_{i=0}^{\infty} z_i$,
where $u_i \in \mathcal{D}(\Y)$ and $v_i \in \mathcal{D}_s(\Y)$ are
the scalar product of rational numbers and connected uni-trivalent
graphs. Let ${\sf D}_{\sf B}$ be the series of all graphs with
components in ${\sf B}$, and ${\sf D}_{\sf C}$ is the series of all
graphs with components in ${\sf C}$. By Lemma \ref{graphs} applied
to ${\sf D}_{\sf B}$ and ${\sf D}_{\sf C}$, we have
$\mathfrak{G}({\sf D}_{\sf B}, {\sf D}_{\sf C}) \longrightarrow
\mathfrak{U} \circledast \mathfrak{G}_{c}({\sf D}_{\sf B}, {\sf
D}_{\sf C})$. Thus, by the Composition Lemma \ref{CompLemma}
$[(\mathfrak{G}({\sf D}_{\sf B},{\sf D}_{\sf C}),\omega)]_{(o;e)} =
[(\mathfrak{U}, \omega_u)]_e \circ [(\mathfrak{G}_{c}({\sf D}_{\sf
B}, {\sf D}_{\sf C}), \omega)]_{(o;e)}$. From \eqref{elsets} we have
$[(\mathfrak{G}({\sf D}_{\sf B}, {\sf D}_{\sf C}), \omega)]_{(o;e)}
= \exp [(\mathfrak{G}_{c}({\sf D}_{\sf B}, {\sf D}_{\sf C}),
\omega)]_{(o;e)}$. It is true that
$[(\mathfrak{G}_{c}(\cdot,\cdot),\omega)]_{(o;e)} =
\Phi\dbbrc{\cdot}{\cdot}$, since the graphs of
$\mathfrak{G}_{c}(\cdot,\cdot)$  are exactly the terms of
$\dbbrc{\cdot}{\cdot}$. From this and Proposition \ref{prop2} we
have $\Phi \dbbr{{\sf D}_{\sf B}}{{\sf D}_{\sf C}} = \exp \Phi
\dbbrc{{\sf D}_{\sf B}}{{\sf D}_{\sf C}}$. Evaluating $u_i =1$ for
all $u_i \in \mathcal{D}(\Y)$ gives $L_{1}^{\mathcal{D}(\Y)} \Phi
\dbbr{{\sf D}_{\sf B}}{{\sf D}_{\sf C}} = \exp
L_{1}^{\mathcal{D}(\Y)} \Phi \dbbrc{{\sf D}_{\sf B}}{{\sf D}_{\sf
C}}$.  But by Proposition \ref{prop1}, this is just
$\bbr{\Phi^{\circ} ({\sf D}_{\sf B})}{\Phi^{\circ} ({\sf D}_{\sf
C})} = \exp \bbrc{ \Phi^{\circ} ({\sf D}_{\sf B})}{\Phi^{\circ}
({\sf D}_{\sf C})}$. So $\bbr{\exp B}{\exp C} = \exp \bbrc{\exp
B}{\exp C}$, giving the result.
\end{proof}

%%%%%%%%%%
%%%%%%%%%%%
%%%%%%%%% IAIN HAS CHANGED BELOW (March 28th)
%%%%%%%%%%

% 8 -----------------------------------------------------------------------------------------------

\section[A Generalization to Diagrammatic Differential Operators]{A Generalization to Diagrammatic Differential Operators}\label{generalization}

So far we have only discussed ``diagrammatic integration''. In this penultimate section we show that our results extend
to the generality of diagrammatic differential operators (\cite{BLT, thurston}). These  are diagrammatic analogues of differential operators and are important in quantum topology.
Perhaps the best known use of diagrammatic differential operators comes from the celebrated wheels and wheeling theorems, first proved in \cite{BLT}. Wheels gives the value of the Kontsevich integral of the unknot $\Omega$ (see Section~\ref{S:EX}), and wheeling states that, in the notation below, $\chi \circ \partial_{\Omega} : \mathcal{B} \rightarrow \mathcal{A} $ is an algebra isomorphism when $\mathcal{B}$ has one colour, the 1-manifold of $\mathcal{A}$ is connected and $\chi$ is the PBW isomorphism of vector-spaces. We will not pursue this further and instead move directly to the combinatorial problem.

The bilinear operator $\bbr{\cdot}{\cdot}$ applied to suitable $g$ and $h$ in
$\mathcal{D}(Y)$ has the property that whenever
 the number of $y$-coloured univalent vertices in $g$ and $h$ do
not match for some $y \in \Y$, then $\bbr{g}{h}=0$. We can relax
this condition and  declare it non-zero only if  all of the
univalent vertices of $g$ are  glued to univalent vertices of $h$.
We extend this bilinearly for all $\mathcal{D}(\Y) \otimes
\mathcal{D}_s(\Y)$ and denote this new operator (which we call a {\em
diagrammatic differential operator})  by $\rbbr{g}{h}$, where
$$
\rbbr{\cdot}{\cdot}: \mathcal{D}(\Y) \otimes \mathcal{D}_s(\Y) \rightarrow \mathcal{D}(\Y).
$$

Similarly, $\rbbrc{\cdot}{\cdot}$ denotes the primitive part of
$\rbbr{\cdot}{\cdot}$. Note that for $u \in \mathcal{D}_s(\Y)$ we
have $\rbbr{1}{u} = u$. The following is an example of this
operator.

%\graph1b
\xy 0;/r.35pc/: (0,0)*{\text{Let } g_1 = }; (12,0)*{ \xy
    (-1,0)*{}="A"; (2,0)*{}="B"; (6,0)*{}="C"; (9,0)*{}="D";
    (-1,0)*{\bullet}; (2,0)*{\bullet}; (6,0)*{\bullet}; (9,0)*{\bullet};
        (0,-2)*{y_1}; (9,-2)*{y_2};
    "A"; "B" **\dir{-};
    "C"; "D" **\dir{-};
    (4,0)*\xycircle(2,2){-}="f";
    "f";
    \endxy};
(24,0)*{ \text{ and } h_1 = }; (36,0)*{ \xy
         (0,3.3)*{y_1};
     (0,-3.3)*{y_1};
    (9,3.3)*{y_2};
    (9,-3.3)*{y_2};
    (0,2)*{}="A"; (2.3,1.3)*{}="B"; (0,-2)*{}="C"; (2.3,-1.3)*{}="D";
    (8,2)*{}="E"; (5.6,1.3)*{}="F"; (8,-2)*{}="G"; (5.6,-1.3)*{}="H";
    (0,2)*{\bullet}; (2.3,1.3)*{\bullet}; (0,-2)*{\bullet}; (2.3,-1.3)*{\bullet};
    (8,2)*{\bullet}; (5.6,1.3)*{\bullet}; (8,-2)*{\bullet}; (5.6,-1.3)*{\bullet};
    "A"; "B" **\dir{-};
    "C"; "D" **\dir{-};
    "E"; "F" **\dir{-};
    "G"; "H" **\dir{-};
    (4,0)*\xycircle(2,2){-}="f";
    "f";
    \endxy};
(43,-1)*{\text{ , then}};
%line 2
(-3,-12)*{(4)}; (10,-12)*{\rbbr{g_1}{h_1} =~}; (16,-12)*{2};
(20,-12)*{\xy
    (0,5)*{}="A"; (4,5)*{}="B";
    (0,-1)*{}="E"; (4,-1)*{}="F";
    (0,5)*{\bullet}; (4,5)*{\bullet};
    (0,-1)*{\bullet}; (4,-1)*{\bullet};
    (-2.5,-4)*{\bullet}; (6.2,-4)*{\bullet};
    (0.3,-3.2)*{\bullet}; (3.7,-3.2)*{\bullet};
    (-2.5,-4)*{}="G"; (6.2,-4)*{}="H";
    (0.3,-3.2)*{}="I"; (3.7,-3.2)*{}="J";
    (2,5)*\xycircle(2,2){-}="f1";
    (2,-2)*\xycircle(2,2){-}="f2";
    "f1"; "f2";
    "A";"E" **\crv{(-2,2)};
    "B";"F" **\crv{(6,2)};
    "I";"G" **\dir{-};
    "J";"H" **\dir{-};
    (-2.5,-5.5)*{y_1}; (6.2,-5.5)*{y_2};
    \endxy};
(27,-12)*{+}; (31,-12)*{2}; (43,-12)*{\xy (-5,2)*{}="A";
(-5,-2)*{}="B";
    (5,2)*{}="C"; (5,-2)*{}="D"; (3,0)*{}="E"; (1,0)*{}="F";
    (7,0)*{}="G"; (9,0)*{}="H";
    (-5,2)*{\bullet}; (-5,-2)*{\bullet};
    (5,2)*{\bullet}; (5,-2)*{\bullet}; (3,0)*{\bullet}; (1,0)*{\bullet};
    (7,0)*{\bullet}; (9,0)*{\bullet};
    (-5,0)*\xycircle(2,2){-}="f1";
    (5,0)*\xycircle(2,2){-}="f2"; "A";"C"**\crv{(0,5)};
    "B";"D"**\crv{(0,-5)}; "E";"F"**\dir{-}; "G";"H"**\dir{-}; (-1,0)*{y_1};
    (11,0)*{y_2};
\endxy
}; (54,-14)*{.};
\endxy

%(-2,0)*{2};
%(0,5)*{}="A"; (4,5)*{}="B";
%(0,-5)*{}="C"; (4,-5)*{}="D";
%(2,5)*\xycircle(2,2){-}="f1";
%(2,0)*\xycircle(2,2){-}="f2";
%(2,-5)*\xycircle(2,2){-}="f3";
%"f1"; "f2";"f3";
%"A";"C" **\dir{-};
%"B";"D" **\dir{-};

When calculating $\rbbr{g}{h}$, we can specify the coloured
univalent vertices of $h$ that will not be identified with the ones
in $g$ by marking them as open vertices ($\circ$) and then
identifying the remaining ones using $\bbr{\cdot}{\cdot}$ (by
definition, $\bbr{\cdot}{\cdot}$ treats open vertices as inert).
After the identification, all the remaining univalent vertices will
be open vertices.  By treating these coloured univalent vertices as
{\em filled vertices} ($\bullet$), we can express $\rbbr{g}{h}$ as a
sum of $\bbr{g}{\cdot}$. We use $a =_{\circ} b$ to indicate that
$a=b$ where the open vertices of $b$ have been filled. For example
from $(4)$,

% {\bf Iain: is $\rbbr{1}{u}=u$ true from your definition? - yes }

%\graph2b
\[
\xy 0;/r.35pc/: (-1,0)*{ \rbbr{g_1}{h_1} =_{\circ}~};
 (20,0)*{ \xy (-1,0)*{\Bigl<};
(6,0)*{ \xy
    (-1,0)*{}="A"; (2,0)*{}="B"; (6,0)*{}="C"; (9,0)*{}="D";
    (-1,0)*{\bullet}; (2,0)*{\bullet}; (6,0)*{\bullet}; (9,0)*{\bullet};
        (0,-2)*{y_1}; (9,-2)*{y_2};
    "A"; "B" **\dir{-};
    "C"; "D" **\dir{-};
    (4,0)*\xycircle(2,2){-}="f";
    "f";
    \endxy};
(14,-1)*{,}; (22,0)*{ \xy
         (0,3.6)*{y_1};
     (0,-3.6)*{y_1};
    (9,3.6)*{y_2};
    (9,-3.6)*{y_2};
    (0,2)*{}="A"; (2.3,1.3)*{}="B"; (0,-2)*{}="C"; (2.3,-1.3)*{}="D";
    (8,2)*{}="E"; (5.6,1.3)*{}="F"; (8,-2)*{}="G"; (5.6,-1.3)*{}="H";
    (0,2)*{\bullet}; (2.3,1.3)*{\bullet}; (2.3,-1.3)*{\bullet};
    (8,2)*{\bullet}; (5.6,1.3)*{\bullet}; (5.6,-1.3)*{\bullet};
    "A"; "B" **\dir{-};
    "E"; "F" **\dir{-};
    "C"; "D" **\dir{-};
    "G"; "H" **\dir{-};
    (-0.6,-2.4)*{\circ};
    (8.6,-2.4)*{\circ};
    (4,0)*\xycircle(2,2){-}="f";
    "f";
    \endxy};
(28,0)*{\Bigr>}; \endxy}; (37,0)*{+}; (55,0)*{ \xy (-1,0)*{\Bigl<};
(6,0)*{ \xy
    (-1,0)*{}="A"; (2,0)*{}="B"; (6,0)*{}="C"; (9,0)*{}="D";
    (-1,0)*{\bullet}; (2,0)*{\bullet}; (6,0)*{\bullet}; (9,0)*{\bullet};
        (0,-2)*{y_1}; (9,-2)*{y_2};
    "A"; "B" **\dir{-};
    "C"; "D" **\dir{-};
    (4,0)*\xycircle(2,2){-}="f";
    "f";
    \endxy};
(14,-1)*{,}; (22,0)*{ \xy
         (0,3.6)*{y_1};
     (0,-3.6)*{y_1};
    (9,3.6)*{y_2};
    (9,-3.6)*{y_2};
    (0,2)*{}="A"; (2.3,1.3)*{}="B"; (0,-2)*{}="C"; (2.3,-1.3)*{}="D";
    (8,2)*{}="E"; (5.6,1.3)*{}="F"; (8,-2)*{}="G"; (5.6,-1.3)*{}="H";
    (2.3,1.3)*{\bullet}; (0,-2)*{\bullet};(2.3,-1.3)*{\bullet};
    (5.6,1.3)*{\bullet}; (8,-2)*{\bullet}; (5.6,-1.3)*{\bullet};
    "A"; "B" **\dir{-};
    "C"; "D" **\dir{-};
    "E"; "F" **\dir{-};
    "G"; "H" **\dir{-};
    (-0.6,2.3)*{\circ};
    (8.6,2.3)*{\circ};
    (4,0)*\xycircle(2,2){-}="f";
    "f";
    \endxy};
(28,0)*{\Bigr>}; \endxy};
%line2
(2,-11)*{+}; (20,-11)*{\xy (-1,0)*{\Bigl<}; (6,0)*{ \xy
    (-1,0)*{}="A"; (2,0)*{}="B"; (6,0)*{}="C"; (9,0)*{}="D";
    (-1,0)*{\bullet}; (2,0)*{\bullet}; (6,0)*{\bullet}; (9,0)*{\bullet};
        (0,-2)*{y_1}; (9,-2)*{y_2};
    "A"; "B" **\dir{-};
    "C"; "D" **\dir{-};
    (4,0)*\xycircle(2,2){-}="f";
    "f";
    \endxy};
(14,-1)*{,}; (22,0)*{ \xy
         (0,3.6)*{y_1};
     (0,-3.6)*{y_1};
    (9,3.6)*{y_2};
    (9,-3.6)*{y_2};
    (0,2)*{}="A"; (2.3,1.3)*{}="B"; (0,-2)*{}="C"; (2.3,-1.3)*{}="D";
    (8,2)*{}="E"; (5.6,1.3)*{}="F"; (8,-2)*{}="G"; (5.6,-1.3)*{}="H";
    (2.3,1.3)*{\bullet}; (2.3,-1.3)*{\bullet};
    (0,2)*{\bullet}; (5.6,1.3)*{\bullet}; (8,-2)*{\bullet}; (5.6,-1.3)*{\bullet};
    "A"; "B" **\dir{-};
    "C"; "D" **\dir{-};
    "E"; "F" **\dir{-};
    "G"; "H" **\dir{-};
    (-0.6,-2.4)*{\circ};
    (8.6,2.3)*{\circ};
    (4,0)*\xycircle(2,2){-}="f";
    "f";
    \endxy};
(28,0)*{\Bigr>}; \endxy}; (37,-11)*{+}; (55,-11)*{\xy
(-1,0)*{\Bigl<}; (6,0)*{ \xy
    (-1,0)*{}="A"; (2,0)*{}="B"; (6,0)*{}="C"; (9,0)*{}="D";
    (-1,0)*{\bullet}; (2,0)*{\bullet}; (6,0)*{\bullet}; (9,0)*{\bullet};
        (0,-2)*{y_1}; (9,-2)*{y_2};
    "A"; "B" **\dir{-};
    "C"; "D" **\dir{-};
    (4,0)*\xycircle(2,2){-}="f";
    "f";
    \endxy};
(14,-1)*{,}; (22,0)*{ \xy
     (0,3.6)*{y_1};
     (0,-3.6)*{y_1};
    (9,3.6)*{y_2};
    (9,-3.6)*{y_2};
    (0,2)*{}="A"; (2.3,1.3)*{}="B"; (0,-2)*{}="C"; (2.3,-1.3)*{}="D";
    (8,2)*{}="E"; (5.6,1.3)*{}="F"; (8,-2)*{}="G"; (5.6,-1.3)*{}="H";
    (2.3,1.3)*{\bullet}; (0,-2)*{\bullet}; (2.3,-1.3)*{\bullet};
    (8,2)*{\bullet}; (5.6,1.3)*{\bullet}; (5.6,-1.3)*{\bullet};
    "A"; "B" **\dir{-};
    "C"; "D" **\dir{-};
    "E"; "F" **\dir{-};
    "G"; "H" **\dir{-};
    (-0.6,2.3)*{\circ};
    (8.6,-2.4)*{\circ};
    (4,0)*\xycircle(2,2){-}="f";
    "f";
    \endxy};
(28,0)*{\Bigr>}; \endxy}; (70,-13)*{.};
\endxy
\]

%{\bf David: As stated this is wrong!  $ \rbbr{g}{h}
%\neq \bbr{g}{\delta(h)}$  as LHS may have univalent vertices while RHS never does.  Instead the identity  $ \delta (\rbbr{g}{h})= \bbr{g}{\delta(h)}$ holds and the argument needs to be modified slightly to include the addition of the "open vertex" legs after the gluing.   Or am I misreading the result?  }

More concisely, for $h \in \mathcal{D}_s(\Y)$ let $\delta(h) \in
\mathcal{D}(\Y)$ be the series of uni-trivalent graphs that can be
obtained from $h$ by opening univalent vertices in all possible
ways. Then $\rbbr{g}{h} =_{\circ} \bbr{g}{\delta(h)}$. Note that
some of the terms of the linear expansion of $\bbr{g}{\delta(h)}$
may be zero. We illustrate this by calculating $\delta(h_1)$ for
$h_1$ in $(4)$,

%\graph2b
\[
\xy 0;/r.30pc/: (-1,0)*{\text{For } h_1=}; (13,0)*{ \xy
     0;/r.35pc/:
         (0,3.3)*{y_1};
     (0,-3.3)*{y_1};
    (9,3.3)*{y_2};
    (9,-3.3)*{y_2};
    (0,2)*{}="A"; (2.3,1.3)*{}="B"; (0,-2)*{}="C"; (2.3,-1.3)*{}="D";
    (8,2)*{}="E"; (5.6,1.3)*{}="F"; (8,-2)*{}="G"; (5.6,-1.3)*{}="H";
    (0,2)*{\bullet}; (2.3,1.3)*{\bullet}; (0,-2)*{\bullet}; (2.3,-1.3)*{\bullet};
    (8,2)*{\bullet}; (5.6,1.3)*{\bullet}; (8,-2)*{\bullet}; (5.6,-1.3)*{\bullet};
    "A"; "B" **\dir{-};
    "C"; "D" **\dir{-};
    "E"; "F" **\dir{-};
    "G"; "H" **\dir{-};
    (4,0)*\xycircle(2,2){-}="f";
    "f";
    \endxy};
(20,-1)*{\text{,}}; (-1,-11)*{\delta(h_1) =_{\circ} }; (10,-11)*{
\xy
        0;/r.35pc/:
         (0,3.3)*{y_1};
     (0,-3.6)*{y_1};
    (9,3.6)*{y_2};
    (9,-3.6)*{y_2};
    (0,2)*{}="A"; (2.3,1.3)*{}="B"; (0,-2)*{}="C"; (2.3,-1.3)*{}="D";
    (8,2)*{}="E"; (5.6,1.3)*{}="F"; (8,-2)*{}="G"; (5.6,-1.3)*{}="H";
    (0,2)*{\bullet}; (2.3,1.3)*{\bullet}; (0,-2)*{\bullet}; (2.3,-1.3)*{\bullet};
    (8,2)*{\bullet}; (5.6,1.3)*{\bullet}; (8,-2)*{\bullet}; (5.6,-1.3)*{\bullet};
    "A"; "B" **\dir{-};
    "C"; "D" **\dir{-};
    "E"; "F" **\dir{-};
    "G"; "H" **\dir{-};
    (4,0)*\xycircle(2,2){-}="f";
    "f";
    \endxy
}; (18,-11)*{+}; (26,-11)*{ \xy
       0;/r.35pc/:
        (0,3.6)*{y_1};
     (0,-3.6)*{y_1};
    (9,3.6)*{y_2};
    (9,-3.6)*{y_2};
    (0,2)*{}="A"; (2.3,1.3)*{}="B"; (0,-2)*{}="C"; (2.3,-1.3)*{}="D";
    (8,2)*{}="E"; (5.6,1.3)*{}="F"; (8,-2)*{}="G"; (5.6,-1.3)*{}="H";
    (2.3,1.3)*{\bullet}; (0,-2)*{\bullet}; (2.3,-1.3)*{\bullet};
    (8,2)*{\bullet}; (5.6,1.3)*{\bullet}; (8,-2)*{\bullet}; (5.6,-1.3)*{\bullet};
    "A"; "B" **\dir{-};
    "C"; "D" **\dir{-};
    "E"; "F" **\dir{-};
    "G"; "H" **\dir{-};
    (4,0)*\xycircle(2,2){-}="f";
    "f";
    (-0.6,2.3)*{\circ};
    \endxy
}; (34,-11)*{+}; (42,-11)*{ \xy
       0;/r.35pc/:
         (0,3.6)*{y_1};
   (0,-3.6)*{y_1};
    (9,3.6)*{y_2};
    (9,-3.6)*{y_2};
    (0,2)*{}="A"; (2.3,1.3)*{}="B"; (0,-2)*{}="C";;(2.3,-1.3)*{}="D";
    (8,2)*{}="E"; (5.6,1.3)*{}="F"; (8,-2)*{}="G"; (5.6,-1.3)*{}="H";
    (0,2)*{\bullet}; (2.3,1.3)*{\bullet}; (2.3,-1.3)*{\bullet};
    (8,2)*{\bullet}; (5.6,1.3)*{\bullet}; (8,-2)*{\bullet}; (5.6,-1.3)*{\bullet};
    "A"; "B" **\dir{-};
    "C"; "D" **\dir{-};
    "E"; "F" **\dir{-};
    "G"; "H" **\dir{-};
    (4,0)*\xycircle(2,2){-}="f";
    "f";
    (-0.6,-2.4)*{\circ};
    \endxy
}; (50,-11)*{+}; (58,-11)*{ \xy
      0;/r.35pc/:
         (0,3.6)*{y_1};
     (0,-3.6)*{y_1};
    (9,3.6)*{y_2};
   (9,-3.6)*{y_2};
    (0,2)*{}="A"; (2.3,1.3)*{}="B"; (0,-2)*{}="C"; (2.3,-1.3)*{}="D";
    (8,2)*{}="E"; (5.6,1.3)*{}="F"; (8,-2)*{}="G"; (5.6,-1.3)*{}="H";
    (0,2)*{\bullet}; (2.3,1.3)*{\bullet}; (0,-2)*{\bullet}; (2.3,-1.3)*{\bullet};
    (8,2)*{\bullet}; (5.6,1.3)*{\bullet}; (5.6,-1.3)*{\bullet};
    "A"; "B" **\dir{-};
    "C"; "D" **\dir{-};
    "E"; "F" **\dir{-};
    "G"; "H" **\dir{-};
    (4,0)*\xycircle(2,2){-}="f";
    "f";
    (8.6,-2.4)*{\circ};
    \endxy
}; (66,-11)*{+}; (74,-11)*{ \xy
       0;/r.35pc/:
         (0,3.6)*{y_1};
     (0,-3.6)*{y_1};
    (9,3.6)*{y_2};
    (9,-3.6)*{y_2};
    (0,2)*{}="A"; (2.3,1.3)*{}="B"; (0,-2)*{}="C"; (2.3,-1.3)*{}="D";
    (8,2)*{}="E"; (5.6,1.3)*{}="F"; (8,-2)*{}="G"; (5.6,-1.3)*{}="H";
    (0,2)*{\bullet}; (2.3,1.3)*{\bullet}; (0,-2)*{\bullet}; (2.3,-1.3)*{\bullet};
    (5.6,1.3)*{\bullet}; (8,-2)*{\bullet}; (5.6,-1.3)*{\bullet};
    "A"; "B" **\dir{-};
    "C"; "D" **\dir{-};
    "E"; "F" **\dir{-};
    "G"; "H" **\dir{-};
    (4,0)*\xycircle(2,2){-}="f";
    "f";
    (8.6,2.3)*{\circ};
    \endxy
};
%line 2
(0,-23)*{+}; (2,-23)*{\Bigl(}; (10,-23)*{ \xy (0,3.6)*{y_1};
     (0,-3.6)*{y_1};
    (9,3.6)*{y_2};
    (9,-3.6)*{y_2};
    (0,2)*{}="A"; (2.3,1.3)*{}="B"; (0,-2)*{}="C"; (2.3,-1.3)*{}="D";
    (8,2)*{}="E"; (5.6,1.3)*{}="F"; (8,-2)*{}="G"; (5.6,-1.3)*{}="H";
    (0,2)*{\bullet}; (2.3,1.3)*{\bullet}; (2.3,-1.3)*{\bullet};
    (8,2)*{\bullet}; (5.6,1.3)*{\bullet}; (5.6,-1.3)*{\bullet};
    "A"; "B" **\dir{-};
    "E"; "F" **\dir{-};
    "C"; "D" **\dir{-};
    "G"; "H" **\dir{-};
    (-0.6,-2.4)*{\circ};
    (8.6,-2.4)*{\circ};
    (4,0)*\xycircle(2,2){-}="f";
    "f";
    \endxy};
(18,-23)*{+}; (26,-23)*{
    \xy
    0;/r.35pc/:
    (0,3.6)*{y_1};
     (0,-3.6)*{y_1};
    (9,3.6)*{y_2};
    (9,-3.6)*{y_2};
    (0,2)*{}="A"; (2.3,1.3)*{}="B"; (0,-2)*{}="C"; (2.3,-1.3)*{}="D";
    (8,2)*{}="E"; (5.6,1.3)*{}="F"; (8,-2)*{}="G"; (5.6,-1.3)*{}="H";
    (2.3,1.3)*{\bullet}; (0,-2)*{\bullet}; (2.3,-1.3)*{\bullet};
    (5.6,1.3)*{\bullet}; (8,-2)*{\bullet}; (5.6,-1.3)*{\bullet};
    "A"; "B" **\dir{-};
    "C"; "D" **\dir{-};
    "E"; "F" **\dir{-};
    "G"; "H" **\dir{-};
    (-0.6,2.3)*{\circ};
    (8.6,2.3)*{\circ};
    (4,0)*\xycircle(2,2){-}="f";
    "f";
    \endxy
}; (34,-23)*{+}; (42,-23)*{
     \xy
     0;/r.35pc/:
     (0,3.6)*{y_1};
     (0,-3.6)*{y_1};
    (9,3.6)*{y_2};
    (9,-3.6)*{y_2};
    (0,2)*{}="A"; (2.3,1.3)*{}="B"; (0,-2)*{}="C"; (2.3,-1.3)*{}="D";
    (8,2)*{}="E"; (5.6,1.3)*{}="F"; (8,-2)*{}="G"; (5.6,-1.3)*{}="H";
    (0,2)*{\bullet}; (2.3,1.3)*{\bullet}; (2.3,-1.3)*{\bullet};
    (5.6,1.3)*{\bullet}; (8,-2)*{\bullet}; (5.6,-1.3)*{\bullet};
    "A"; "B" **\dir{-};
    "C"; "D" **\dir{-};
    "E"; "F" **\dir{-};
    "G"; "H" **\dir{-};
    (-0.6,-2.4)*{\circ};
    (8.6,2.3)*{\circ};
    (4,0)*\xycircle(2,2){-}="f";
    "f";
    \endxy
}; (50,-23)*{+}; (58,-23)*{
     \xy
     0;/r.35pc/:
         (0,3.6)*{y_1};
     (0,-3.6)*{y_1};
    (9,3.6)*{y_2};
    (9,-3.6)*{y_2};
    (0,2)*{}="A"; (2.3,1.3)*{}="B"; (0,-2)*{}="C"; (2.3,-1.3)*{}="D";
    (8,2)*{}="E"; (5.6,1.3)*{}="F"; (8,-2)*{}="G"; (5.6,-1.3)*{}="H";
    (2.3,1.3)*{\bullet}; (0,-2)*{\bullet};(2.3,-1.3)*{\bullet};
    (8,2)*{\bullet}; (5.6,1.3)*{\bullet}; (5.6,-1.3)*{\bullet};
    "A"; "B" **\dir{-};
    "C"; "D" **\dir{-};
    "E"; "F" **\dir{-};
    "G"; "H" **\dir{-};
    (-0.6,2.3)*{\circ};
    (8.6,-2.4)*{\circ};
    (4,0)*\xycircle(2,2){-}="f";
    "f";
    \endxy
}; (66,-23)*{\Bigr)}; (68,-23)*{+}; (76,-23)*{ \xy
      0;/r.35pc/:
         (0,3.6)*{y_1};
     (0,-3.6)*{y_1};
    (9,3.6)*{y_2};
    (9,-3.6)*{y_2};
    (0,2)*{}="A"; (2.3,1.3)*{}="B"; (0,-2)*{}="C"; (2.3,-1.3)*{}="D";
    (8,2)*{}="E"; (5.6,1.3)*{}="F"; (8,-2)*{}="G"; (5.6,-1.3)*{}="H";
    (0,2)*{\bullet}; (2.3,1.3)*{\bullet}; (0,-2)*{\bullet}; (2.3,-1.3)*{\bullet};
    (5.6,1.3)*{\bullet}; (5.6,-1.3)*{\bullet};
    "A"; "B" **\dir{-};
    "C"; "D" **\dir{-};
    "E"; "F" **\dir{-};
    "G"; "H" **\dir{-};
    (4,0)*\xycircle(2,2){-}="f";
    "f";
    (8.6,2.3)*{\circ};
    (8.6,-2.4)*{\circ};
    \endxy
}; (84,-23)*{+}; (92,-23)*{ \xy
       0;/r.35pc/:
    (0,3.6)*{y_1};
    (0,-3.6)*{y_1};
    (9,3.6)*{y_2};
    (9,-3.6)*{y_2};
    (0,2)*{}="A"; (2.3,1.3)*{}="B"; (0,-2)*{}="C"; (2.3,-1.3)*{}="D";
    (8,2)*{}="E"; (5.6,1.3)*{}="F"; (8,-2)*{}="G"; (5.6,-1.3)*{}="H";
    (2.3,1.3)*{\bullet}; (2.3,-1.3)*{\bullet};
    (8,2)*{\bullet}; (5.6,1.3)*{\bullet}; (8,-2)*{\bullet}; (5.6,-1.3)*{\bullet};
    "A"; "B" **\dir{-};
    "C"; "D" **\dir{-};
    "E"; "F" **\dir{-};
    "G"; "H" **\dir{-};
    (4,0)*\xycircle(2,2){-}="f";
    "f";
    (-0.6,2.3)*{\circ};
    (-0.6,-2.4)*{\circ};
    \endxy
};
%line 3
(0,-35)*{+}; (10,-35)*{ \xy
        0;/r.35pc/:
         (0,3.6)*{y_1};
    (0,-3.6)*{y_1};
    (9,3.6)*{y_2};
    (9,-3.6)*{y_2};
    (0,2)*{}="A"; (2.3,1.3)*{}="B"; (0,-2)*{}="C"; (2.3,-1.3)*{}="D";
    (8,2)*{}="E"; (5.6,1.3)*{}="F"; (8,-2)*{}="G"; (5.6,-1.3)*{}="H";
    (0,2)*{\bullet}; (2.3,1.3)*{\bullet}; (2.3,-1.3)*{\bullet};
    (5.6,1.3)*{\bullet}; (5.6,-1.3)*{\bullet};
    "A"; "B" **\dir{-};
    "C"; "D" **\dir{-};
    "E"; "F" **\dir{-};
    "G"; "H" **\dir{-};
    (4,0)*\xycircle(2,2){-}="f";
    "f";
    (8.6,2.3)*{\circ};
    (8.6,-2.4)*{\circ};
    (-0.6,-2.4)*{\circ};
    \endxy
}; (18,-35)*{+}; (26,-35)*{ \xy
       0;/r.35pc/:
        (0,3.6)*{y_1};
     (0,-3.6)*{y_1};
    (9,3.6)*{y_2};
    (9,-3.6)*{y_2};
    (0,2)*{}="A"; (2.3,1.3)*{}="B"; (0,-2)*{}="C"; (2.3,-1.3)*{}="D";
    (8,2)*{}="E"; (5.6,1.3)*{}="F"; (8,-2)*{}="G"; (5.6,-1.3)*{}="H";
    (2.3,1.3)*{\bullet}; (0,-2)*{\bullet}; (2.3,-1.3)*{\bullet};
    (5.6,1.3)*{\bullet}; (5.6,-1.3)*{\bullet};
    "A"; "B" **\dir{-};
    "C"; "D" **\dir{-};
    "E"; "F" **\dir{-};
    "G"; "H" **\dir{-};
    (4,0)*\xycircle(2,2){-}="f";
    "f";
     (-0.6,2.3)*{\circ};
    (8.6,2.3)*{\circ};
    (8.6,-2.4)*{\circ};
    \endxy
}; (34,-35)*{+}; (42,-35)*{ \xy
       0;/r.35pc/:
       (0,3.6)*{y_1};
    (0,-3.6)*{y_1};
    (9,3.6)*{y_2};
    (9,-3.6)*{y_2};
    (0,2)*{}="A"; (2.3,1.3)*{}="B"; (0,-2)*{}="C"; (2.3,-1.3)*{}="D";
    (8,2)*{}="E"; (5.6,1.3)*{}="F"; (8,-2)*{}="G"; (5.6,-1.3)*{}="H";
    (2.3,1.3)*{\bullet}; (2.3,-1.3)*{\bullet};
    (8,-2)*{\bullet}; (5.6,1.3)*{\bullet}; (5.6,-1.3)*{\bullet};
    "A"; "B" **\dir{-};
    "C"; "D" **\dir{-};
    "E"; "F" **\dir{-};
    "G"; "H" **\dir{-};
    (4,0)*\xycircle(2,2){-}="f";
    "f";
    (-0.6,2.3)*{\circ};
    (-0.6,-2.4)*{\circ};
    (8.6,2.4)*{\circ};
    \endxy
}; (50,-35)*{+}; (58,-35)*{ \xy
      0;/r.35pc/:
        (0,3.6)*{y_1};
     (0,-3.6)*{y_1};
    (9,3.6)*{y_2};
    (9,-3.6)*{y_2};
    (0,2)*{}="A"; (2.3,1.3)*{}="B"; (0,-2)*{}="C"; (2.3,-1.3)*{}="D";
    (8,2)*{}="E"; (5.6,1.3)*{}="F"; (8,-2)*{}="G"; (5.6,-1.3)*{}="H";
    (8,2)*{\bullet}; (2.3,1.3)*{\bullet}; (2.3,-1.3)*{\bullet};
    (5.6,1.3)*{\bullet}; (5.6,-1.3)*{\bullet};
    "A"; "B" **\dir{-};
    "C"; "D" **\dir{-};
    "E"; "F" **\dir{-};
    "G"; "H" **\dir{-};
    (4,0)*\xycircle(2,2){-}="f";
    "f";
    (-0.6,2.3)*{\circ};
    (-0.6,-2.4)*{\circ};
    (8.6,-2.4)*{\circ};
    \endxy
}; (66,-35)*{+}; (74,-35)*{ \xy
       0;/r.35pc/:
    (0,3.6)*{y_1};
     (0,-3.6)*{y_1};
    (9,3.6)*{y_2};
    (9,-3.6)*{y_2};
    (0,2)*{}="A"; (2.3,1.3)*{}="B"; (0,-2)*{}="C"; (2.3,-1.3)*{}="D";
    (8,2)*{}="E"; (5.6,1.3)*{}="F"; (8,-2)*{}="G"; (5.6,-1.3)*{}="H";
    (2.3,1.3)*{\bullet}; (2.3,-1.3)*{\bullet};
    (5.6,1.3)*{\bullet}; (5.6,-1.3)*{\bullet};
    "A"; "B" **\dir{-};
    "C"; "D" **\dir{-};
    "E"; "F" **\dir{-};
    "G"; "H" **\dir{-};
    (4,0)*\xycircle(2,2){-}="f";
    "f";
    (-0.6,2.3)*{\circ};
    (-0.6,-2.4)*{\circ};
    (8.6,-2.4)*{\circ};
    (8.6,2.3)*{\circ};
    \endxy
}; (81,-37)*{.};
\endxy
\]

Note that the terms in parenthesis are the ones that give a nonzero
contribution in  $\bbr{g_1}{\delta(h_1)}$. Using the above
observation we can prove an analogue of the main result.

\begin{corollary} \label{cor2}
Let $B,C\in \mathcal{D}_s(\Y)$ be strutless and primitive, then
\begin{equation*}
\rbbr{\exp B}{\exp C} = \exp\left( \rbbrc{\exp B}{\exp C} \right).
\end{equation*}
\end{corollary}

\begin{proof}

Notice that if $C \in \mathcal{D}_s(\Y)$ is strutless and primitive,
so is $\delta(C)$. It is clear that $\rbbr{\exp B}{\exp C} =_{\circ}
\bbr{\exp B}{\delta (\exp C)}= \bbr{\exp B}{\exp \delta (C)}$. But
by Theorem \ref{thm2} we have the relation $\bbr{\exp B}{\exp \delta
(C)} = \exp \bbrc{\exp B}{\exp \delta (C)}$. But $\exp \rbbrc{\exp
B}{\exp C} =_{\circ} \exp \bbrc{\exp B}{\exp \delta (C)} $,
completing the proof.
\end{proof}

We hope that Corollary~\ref{cor2} may prove useful for finding
expressions for values of the Kontsevich invariant in algebras other
than $\mathcal{B}$, perhaps through the use of the wheeling theorem
mentioned above.

% 9 -----------------------------------------------------------------------------------------------

\section{Examples}\label{Examples}\label{S:EX}

Although the LMO invariant can be computed algorithmically to any finite degree, there are  few known examples of the full values of this invariant.  Known explicit examples include lens spaces (\cite{LMO, BNR}) and certain Seifert fiber spaces (\cite{BNR}). As some applications of our results we shall use Theorem~\ref{thm2}, its corollaries and results of Bar-Natan and Lawrence to determine the logarithm of the LMO invariant of certain manifolds. The logarithm of the LMO invariant is known as the  {\em primitive LMO invariant}, and is denoted by $z^{\text{LMO}}$.

The principle advantage of looking at primitive finite-type and quantum invariants is that their structure and the coefficients of their terms are often more accessible than the original invariant (\cite{ohtsuki}). Therefore  primitive invariants and the corresponding space of primitive diagrams are well studied in knot theory.
In addition to this the primitive LMO invariant is known to  behave well under the connect sum operations of 3-manifolds and  reversal of orientation (\cite{LMO}). For example, if $M$ and $N$ are two rational homology spheres and $M \# N$ their connected sum,   then $z^{LMO} (M \# N) =  z^{LMO} (M ) + z^{LMO} (N)$ (to see this note that the framed links representing $M$ and $N$ have disjoint colouring so the formal Gaussian integration can be carried out for each set of  variables  separately). This formula can be applied to the formulae below to obtain expressions for the sums of the manifolds, although we do not include details here.

Through clever use of the wheels and wheeling formulae, Bar-Natan and Lawrence, in \cite{BNR}, gave explicit calculations of the Kontsevich integral of integrally framed Hopf links and Hopf chains. Using these calculations  they  went on to calculate the LMO invariant of lens spaces, which may be presented as integrally framed Hopf chains (\cite{Rolfsen}) and certain Seifert fiber spaces which have a simple ``key chain'' presentation (\cite{Mo, scott}).  By considering these results, we use our formulae to calculate the primitive LMO invariants of these manifolds.

Let $\omega_{2n}$ be the {\em wheel} of degree $2n$, {\em ie.} the
uni-trivalent graph  made from a $2n$-gon with an additional edge
coming out from each vertex. We assume $\omega$ has $x$-coloured
univalent vertices. Also let $\Omega_x = \exp \left(
\sum_{m=1}^{\infty} b_{2m} \; \omega_{2m} \right)$
 denote the Kontsevich integral of the unknot  and
 $\Omega_{x/p} =  \exp \left(  \sum_{m=1}^{\infty} b_{2m}/p^{2m} \;
\omega_{2m} \right)$,
  where  the $b_{2m} \in \mathbb{Q}$
 are the modified Bernoulli numbers (see \emph{e.g.} \cite{BLT}). Finally $\theta$
denotes the planar trivalent graph with two vertices. All vertex orientations are inherited from the plane.

Bar-Natan and Lawrence show (in the proof of their Proposition~5.1) that the LMO invariant of
the~$(p,q)$ lens space is given by the formula
\[
Z^{\text{LMO}}(L_{p,q}) =  \exp \left(  \frac{-S(q/p)}{48} \; \theta   \right) \bbr{\Omega_x}{\Omega_x}^{-1} \bbr{\Omega_x}{\Omega_{x/p}},
\]
where $S(q/p) \in \mathbb{Q}$ is the Dedekind symbol (\cite{KM}), whose definition we do not include here.
We may either apply Theorem~\ref{thm2} to this formula and obtain
\[
Z^{\text{LMO}}(L_{p,q}) = \exp \left(  \frac{-S(q/p)}{48} \; \theta   \right)
\exp \left(  \bbrc{\Omega_x}{\Omega_x}  \right)^{-1}
\exp \left(  \bbrc{\Omega_x}{\Omega_{x/p}}  \right),
\]
or one can use Corollary~3.5 of \cite{BNR} to write
$ \bbr{\Omega_x}{\Omega_x}^{-1} \bbr{\Omega_x}{\Omega_{x/p}}$ as
$\bbr{\Omega_x}{\Omega_x^{-1}  \Omega_{x/p}}$, then apply Theorem~\ref{thm2} to get
\[
Z^{\text{LMO}}(L_{p,q}) = \exp \left(  \frac{-S(q/p)}{48}\; \theta   \right)
\exp \left( \bbrc{\Omega_x}{\Omega_x^{-1}  \Omega_{x/p}} \right).
\]
Since $\mathcal{B}$ is a commutative algebra we obtain the following.
\begin{proposition}
The primitive LMO invariant of a (p,q) lens space is given by
\begin{eqnarray*}
z^{\text{LMO}}(L_{p,q}) =
\bbrc{\Omega_x}{\Omega_{x/p}}
-\bbrc{\Omega_x}{\Omega_x}
-\frac{S(q/p)}{48} \; \theta
\\
=
\bbrc{\Omega_x}{\Omega_x^{-1}  \Omega_{x/p}} - \frac{S(q/p)}{48} \; \theta .
\end{eqnarray*}
\end{proposition}

As our concluding example, if $M = S^3(b, p_1/q_1, \ldots , p_n/q_n )$ is the
Seifert fibered space with a spherical base described in Section~5.2 of \cite{BNR}, one can use Bar-Natan and Lawrence's formula and Theorem~\ref{thm2} in a similar way to calculate its primitive LMO invariant as
\begin{eqnarray*}
z^{\text{LMO}}(M) =
\left\langle \exp  \left( \frac{1}{2e_0} \, \cstrut{x}{x} \right) \,
, \; \Omega_{x}^{2-n}
\prod_{i} \Omega_{x/p_i} \right\rangle_c
 -\langle \Omega_x ,  \Omega_{x} \rangle_c
\\
+  \frac{1}{4} \left(
\lambda_{\omega}(M)
+\frac{1}{12 e_0}
 \left(  n-2- \sum_{i} \frac{1}{p_i^2}  \right)
\right) \theta ,
\end{eqnarray*}
where
$\lambda_{\omega}(M)$ denotes the Casson-Walker invariant (\cite{walker})   of $M$ and
$e_0 := b + \sum_i q_i /p_i$.

\medskip
}   %%%%%%%%%%%%%%%%%%% parskip at beginning %%%%%%%%%%%%%%%%%%%%%%%%%%%%%%%%%}
%%%%%%%%%%%%%%%%%%%%%%%%%%%%%%%%%%%%%%%%%%%%%%%%%%%%%%%%%%%%%%%%%%%%%%%%%%

\end{document}